\documentclass[a4paper,11pt]{article}

\newtheorem{thm}{Theorem}[section]

\newtheorem{dfn}[thm]{Definition}

\usepackage[all]{xy}

\newcommand{\blue}{\bf }

\newcommand{\co}{\colon}

  \newcommand{\clpse}{\mathsf{clpse}}
\usepackage{amssymb}
\usepackage{amsbsy}
\usepackage{amsmath}

\newcommand{\ON}{\ensuremath{\mathrm{ON}}}  

\renewcommand{\phi}{\varphi}

\newcommand{\nstrong}{\ensuremath{\not\kern-4pt\lhd\;}} 

\newbox\noforkbox \newdimen\forklinewidth
\forklinewidth=0.3pt \setbox0\hbox{$\textstyle\smile$}
\setbox1\hbox to \wd0{\hfil\vrule width \forklinewidth depth-2pt
  height 10pt \hfil}
\wd1=0 cm \setbox\noforkbox\hbox{\lower 2pt\box1\lower
2pt\box0\relax}
\def\unionstick{\mathop{\copy\noforkbox}\limits}

\def\nonfork_#1{\unionstick_{\textstyle #1}}

\setbox0\hbox{$\textstyle\smile$} \setbox1\hbox to \wd0{\hfil{\sl
/\/}\hfil} \setbox2\hbox to \wd0{\hfil\vrule height 10pt depth
-2pt width
                \forklinewidth\hfil}
\newbox\doesforkbox
\setbox\doesforkbox\hbox{\lower 2pt\box1 \lower
2pt\box2\lower2pt\box0\relax}
\def\nunionstick{\mathop{\copy\doesforkbox}\limits}

\def\fork_#1{\nunionstick_{\textstyle #1}}

\newcommand{\comment}[1]{}  

\newtheorem{cor}[thm]{Corollary}

\newtheorem{lem}[thm]{Lemma}
\newtheorem{deff}[thm]{Definition}

\newtheorem{conj}[thm]{Conjecture}

\newcommand{\ZFC}{{\mathbf{ZFC}}}

\newcommand{\PA}{{\mathbf{PA}}}

\newcommand{\prf}{{\bf Proof: }}

\newcommand{\gdw}{\,\leftrightarrow\,}

\newcommand{\BI}{{\mathbf{BI}}}

\newcommand{\beq}{\begin{eqnarray}}
\newcommand{\eeq}{\end{eqnarray}}

\newcommand{\paar}[1]{\langle #1\rangle}

\newcommand{\ran}{{\mathrm{ran}}}

\newcommand{\Ord}{{\mathbb{ORD}}}

\newcommand{\PAs}[1]
{#1^+}

\newcommand{\bes}{\begin{eqnarray*}}
\newcommand{\ees}{\end{eqnarray*}}

\newcommand{\PRA}{\mathbf{PRA}}

\newcommand{\ACA}{{\mathbf{ACA}}}

\newcommand{\DI}{\displaystyle}

\newcommand{\OR}{{\cal{OR}}}

\newcommand{\Cut}{$(Cut)$}

\newcommand{\WO}{\mathrm{WO}}

\newcommand{\TI}{{\mathbf {TI}}}

\newcommand{\CNF}{=_{_{CNF}}}
 
\newcommand{\OM}{\Omega}
\newcommand{\Bach}{{\mathfrak B}}
\newcommand{\fiedi}{\varphi^{\!^{\mathfrak B}}}
\newcommand{\Zz}{{\mathbf{Z_2}}}

\newcommand{\fX}{{\mathfrak X}}

\newcommand{\WOP}{{\mathbf{WOP}}}
\newcommand{\WOPP}{{\mathbb{HWOP}}}
\newcommand{\WOPPP}{{\mathbb{HWOP}}_1}
\newcommand{\TUQ}{{\mathcal T}^{\infty}_Q}
\newcommand{\fXO}{{\mathfrak X}_0}

\newcommand{\LOO}{{\mathbb{LO}}}

\newcommand{\kommentar}[1]{}

\newcommand{\LIO}{\mathbb{LO}}
\newcommand{\supp}{\mathsf{supp}}
\newcommand{\EEO}{\mathsf{supp}_{\Omega}}
\newcommand{\EEON}{\mathsf{supp}_{\Omega_n}}
\newcommand{\EEOX}{\mathsf{supp}_{\Omega}^{\!_{\mathfrak X}}}
\newcommand{\BHS}{\mathsf{OT}(\vartheta)}
\newcommand{\BHSX}{\mathsf{OT}_{\!_{\mathfrak X}}(\vartheta)}
\newcommand{\varthetao}{\vartheta}
\newcommand{\varthetax}{\vartheta_{\!_{\mathfrak X}}}
\newcommand{\varthetad}{\vartheta_{\!_{\mathcal D}}}
\newcommand{\Eb}{\mathfrak{E}}
\newcommand{\vieh}{\varphi}
\newcommand{\EDD}{\mathfrak{E}}
\newcommand{\BHSD}{\mathsf{OT}_{\!_{\mathcal D}}(\vartheta)}
\newcommand{\EEOD}{\mathsf{supp}_{\Omega}^{\!_{\mathcal D}}}
\newcommand{\PAO}{\mathbf{PA}_{_{\mathfrak X}}}
\newcommand{\PAI}{\mathbf{PA}^{\!*}_{_{\mathfrak I}}}
\newcommand{\PAIs}{\mathbf{PA}_{_{\mathfrak I}}}
\newcommand{\LQ}{L_{_\mathfrak{X}}}
\newcommand{\DeltaO}{\Delta_{_\Omega}}
\newcommand{\SigmaO}{\Sigma_{_\Omega}}
\newcommand{\bbM}[1]{{\mathbb #1}}

\newcommand{\TsQ}[3]{\provx{{\mathbb T}_{\mathfrak I}^{^{\bar{Q}}}}{#1}{#2}{#3}}

\usepackage{latexsym}
\usepackage{prooft}
\usepackage{bussproofs}

\begin{document}

\title{Well-Ordering Principles in Proof Theory and  Reverse Mathematics} 

\author{Michael Rathjen\\
{\small Department of
Pure Mathematics, University of Leeds}\\ {\small Leeds LS2 9JT,
UK}\\ {\small E-Mail: M.Rathjen@leeds.ac.uk}}

\date{}
\maketitle 

\begin{abstract} 
Several theorems about the equivalence of familiar theories of reverse mathematics
with certain well-ordering principles have been proved
by recursion-theoretic and combinatorial methods  (Friedman, Marcone, Montalb\'an et al.) and with far-reaching results 
 by proof-theoretic  technology
(Afshari, Freund, Girard, Rathjen, Thomson,  Valencia Vizca\'ino, Weiermann et al.), employing deduction search trees and cut elimination 
theorems in infinitary logics with ordinal bounds in the latter case. 

At type level one, the well-ordering principles are of the form 
$$(*)\mbox{  "if $X$ is well-ordered then $f(X)$ is well-ordered"}$$  where $f$ is a standard proof theoretic 
function from ordinals to ordinals (such $f$'s are always dilators).
One aim of the paper is to present a general methodology underlying these results that enables one to construct omega-models of particular theories from $(*)$ and even $\beta$-models 
from the type 2 version of $(*)$.

As $(*)$  is of complexity $\Pi^1_2$ such a principle cannot characterize stronger comprehensions at the level of $\Pi^1_1$-comprehension. This requires a higher order version of $(*)$ that employs ideas from ordinal representation systems with collapsing functions used in impredicative proof theory.  The simplest one is the Bachmann construction. Relativizing the latter construction to any dilator $f$ and asserting that this always yields a well-ordering turns out to be equivalent to $\Pi^1_1$-comprehension. This result has been conjectured and a proof strategy adumbrated roughly 10 years ago, but the proof has only been worked out by Anton Freund in recent years.
\comment{\\[2ex]
Several results about the equivalence of familiar theories of reverse mathematics
with certain well-ordering principles have been proved (Friedman, Marcone, Montalb\'an et al.)
by recursion-theoretic and combinatorial methods and also by proof theory
(Afshari, Girard, Rathjen, Thomson,  Valencia Vizca\'ino, Weiermann et al.), employing deduction search trees and cut elimination 
theorems in infinitary logics with ordinal bounds.
 
One goal of this survey paper is to present a general methodology underlying these results which in many cases allows one to establish an equivalence between two types of statements. The first type is concerned with the existence of omega models of a theory whereas the second type asserts that a certain (usually well-known) elementary operation on orderings preserves the property of being well-ordered.
These operations are related to ordinal representation systems (ORS) that play a central role in proof theory. The question of naturality of ORS has vexed logicians for a long time. While ORS have a low computational complexity their ``true" nature evades characterization in those terms. One attempt has been to describe their structural properties in category-theoretic terms (Aczel, Feferman, Girard et al.).
These ideas will be discussed in this article. 

A second goal is to present rather recent developments (due to Arai, Freund, Rathjen), especially work by Freund on higher order well-ordering principles and comprehension.}
\\[2ex]
Keywords:  infinite proof theory, reverse mathematics, ordinal analysis, ordinal representation systems, cut elimination, dilators, $\omega$-models, $\beta$-models
 \\ 
 MSC2000: 03F50;  03F25;
03E55;  03B15; 03C70
\end{abstract}

\section{Introduction} The aim to illuminate the role of the infinite in mathematics gave rise to set theory and proof theory alike.
Whereas  set theory takes much of the infinite for granted (e.g. full separation and powerset), proof theory 
strives to analyze it from a stance of potential infinity. While the objects that proof theory primarily deals with are rather concrete (e.g. theories, proofs, ordinal representation systems),  it is also concerned with abstract higher type\footnote{Or {\em ideal} properties in Hilbert's sense.} properties of concrete objects (e.g. well-foundedness, preservation of well-foundedness).
 In ordinal analysis the impetus is to associate an ordinal representation system $\OR$ with a theory $T$ in such a way that the former displays the commitments to the infinite encapsulated in $T$. Less poetically,  it entails that $\PRA+\TI_{qf}(\OR)$ proves all $\Pi^0_2$-theorems of $T$, however, no proper initial segment of $\OR$ suffices for that task (where $\PRA$ refers to primitive recursive arithmetic and $\TI_{qf}(\OR)$ stands for quantifier-free transfinite induction along $\OR$).\footnote{In some books, the proof-theoretic ordinal $|T|$ of a theory $T$ is defined as the supremum of the order-types ${\mid}\prec{\mid}$ of the $T$-provable  well-orderings $\prec$, i.e., 
 $$|T|_{\sup} \;=\; \sup \{|{\prec}| \,:\, \mbox{$\prec$ primitive recursive }; T\vdash \mathrm{WO}(\prec)\}$$
 where $\mathrm{WO}$ expresses that $\prec$ is a well-ordering. This definition, however, conflates extensional aspects, i.e. the order-type of $\prec$ in some background universe of set theory, with intensional aspects, i.e.
 the niceties of the coding of $\prec$ in the integers. This is liable to lead to results where $|T|_{\sup}$ has nothing to do with the ``real" proof-theoretic strength of $T$. For example, one can
 define ``faithful" theories $\mathsf{T}_1$ and $\mathsf{T}_2$ such that:
(a) $|\mathsf{T}_1|_{\sup}$ is the Bachmann-Howard ordinal but $T_1$ is finitistically reducible to $\PA$, proving the same $\Pi^0_1$-theorems as $\PA$;
(b) $|\mathsf{T}_2|_{\sup} = \varepsilon_0$ but $T_2$ proves the same arithmetical theorems as an extension of $\ZFC$ with large cardinal axioms.
Thus, paraphrasing the title of Paul Benacerraf's famous 1965 paper, one could say that $|T|_{\sup}$ displays what the proof-theoretic ordinal of $T$ could not be.
For more on this, see \cite{Rathjen Realm}.} The first ordinal representations were offshoots of ordinal normal forms such as Cantor's and Veblen's from more than a hundred years ago. Gentzen employed the Cantor normal form with base $\omega$ to provide an ordinal representation system for $\varepsilon_0$ in his last two papers \cite{gentzen38,gentzen43}, which mark the beginning of ordinal-theoretic proof theory.

It is often stressed that ordinal representation systems are computable\footnote{They actually possess a very low computational complexity; see \cite{so1}.}  structures, which is true and allows for the treatment of their order-theoretic and algebraic aspects in very weak systems of arithmetic, but it is only one of their distinguishing features. Overstating the computability aspect tends to give the impression that their study is part of the venerable research area of computable orderings. In actual fact, the two subjects have very little in common. 

 There are many specific ordinal representations systems to be found in proof theory. To make the choice of these systems intelligible one would like to discern general principles 
involved. A particular challenge is posed by making their algebraic features explicit, i.e. those features not related to effectiveness or well-foundedness. This constitutes a crucial step towards a general theory via an axiomatic approach.
As far as I'm aware, it was Feferman \cite{f68}  
 who first initiated such a theory by isolating the property 
of {\em effective relative categoricity} as the fulcrum of such a characterization. Through the notion of relative categoricity he succeeded in crystallizing the algebraic aspect of ordinal representation systems by way of relativizing them to any set of order-indiscernibles, thereby in effect scaling them up to functors on the category of linear orders with order-preserving maps as morphisms. 

Ordinal representation systems understood in this relativized way give rise to well-ordering principles: $$\mbox{``if $\mathfrak X$ is a well-ordering then so is $\OR_{\mathfrak X}$"}$$ where
 $\OR_{\mathfrak X}$ arises from $\OR$ by letting $\mathfrak X$ play the role of the order-indiscernibles. Several theories of reverse mathematics have been shown to be equivalent to
such principles involving iconic ordinal representation system from the proof-theoretic literature.
From a technical point of view the methods can be roughly divided into  recursion-theoretic and combinatorial tools  (Friedman, Marcone, Montalban et al.) on the one hand and 
proof-theoretic tools 
(Afshari, Freund, Girard, Rathjen, Thomson,  Valencia Vizca\'ino, Weiermann et al.) on the other hand, where in the latter approach the employment of 
 deduction search trees and cut elimination 
theorems in infinitary logics with ordinal bounds play a central role. 

This is a survey paper about these results that tries to illuminate the underlying ideas, emphasizing the proof theory side.

\subsection{Reverse mathematics} It is assumed that the reader is roughly familiar the the program of {\em Reverse Mathematics} (RM) and the language and axioms of the formal system
of second order arithmetic, $\Zz$ as for instance presented in the standard reference \cite{SOSA}.
Just as a reminder, RM started out with the observation that for many mathematical theorems $\tau$, there is a weakest natural
subsystem $S(\tau)$ of $\Zz$
 such that $S(\tau)$ proves $\tau$.
  Moreover, it has turned out that $S(\tau)$
often belongs to a small list of specific subsystems of $\Zz$.
In particular, Reverse Mathematics  has singled out five subsystems of
$\Zz$, often referred to as the {\em big five}, that provide (part of) a standard scale against which the strengths of theories can often be measured: 
\begin{enumerate} 
\item  ${\mathbf{RCA}}_0$ \phantom{AAAAAA}  Recursive
Comprehension
\item  ${\mathbf{WKL}}_0$ \phantom{AAAAAA} Weak K\"onig's
Lemma
 \item   ${\mathbf{ACA}}_0$ \phantom{AAAAAA}
Arithmetical Comprehension
  \item   ${\mathbf{ATR}}_0$ \phantom{AAAAAA}
Arithmetical Transfinite Recursion
 \item  $\Pi^1_1{-}{\mathbf{CA}}_0$ \phantom{AAA}
 $\Pi_1^1$-Comprehension
 \end{enumerate}

\section{History}
The principles we  will be concerned with are  particular $\Pi^1_2$ statements
of the form
 \begin{eqnarray}\label{1}\WOP(f):&&\;\,\,\;\;\;\forall X\,[\WO(\fX)\rightarrow \WO(f(\fX))]\end{eqnarray} where $f$ is a standard proof-theoretic function from ordinals to ordinals and $\WO(\fX)$ stands for `$\fX$ is a well-ordering'.\footnote{And $\WOP$ is an acronym for the German word ``{\bf W}ohl{\bf o}rdnungs{\bf p}rinzip''.}
 There are by now several examples of functions $f$ familiar from proof theory where the statement
 $\WOP(f)$ has turned out to be equivalent to one of the theories of reverse
 mathematics over a weak base theory (usually ${\mathbf{RCA}}_0$). The first explicit example appears to be due to Girard
 \cite[5.4.1 theorem] {girard87} (see also \cite{Hirst}). However, it is also implicit in Sch\"utte's
 proof of cut elimination for $\omega$-logic \cite{sch51} and ultimately has its roots in Gentzen's work, namely in his first unpublished consistency proof,\footnote{The original German version was finally published in 1974 \cite{gentzen74}. An earlier English
translation appeared in 1969 \cite{gentzen69}.}  where he introduced the notion of a ``Reduziervorschrift" \cite[p. 102]{gentzen74} for a sequent. The latter is
a well-founded tree built bottom-up via ``Reduktionsschritte", starting with the given sequent
and passing up from conclusions to premises until an axiom is reached.

\subsection{ $2^{\mathfrak X}$ and Arithmetical Comprehension}
 In our first example, one constructs a linear ordering $2^{\mathfrak X} = ( |2^X|,<_{_{2^X}})$ from a given linear ordering ${\mathfrak X}=(X,<_X)$, where $|2^X|$ consists of the formal
 sums $2^{x_1}+\ldots+ 2^{x_n}$ with $x_n<_{\mathfrak X} \ldots <_{\mathfrak X}x_1$, and the ordering between  two formal sums $2^{x_0}+\ldots+ 2^{x_n}$, $2^{y_0}+\ldots+ 2^{y_m}$ is determined 
 as follows: \begin{eqnarray*} 2^{x_0}+\ldots+ 2^{x_n} <_{2^X}2^{y_1}+\ldots+ 2^{y_m} &\mbox{ iff }& \exists i\leq \max(m,n)[x_i<_X y_i\,\wedge\, \forall j<i\;x_j=y_j]\\
     &&\mbox{ or } n<m \,\wedge\, \forall i\leq n \; x_i=y_i.\end{eqnarray*}

 \begin{thm}\label{girard} {\em (Girard 1987)}
Over ${\mathbf{RCA}}_0$ the following are equivalent:
 \begin{itemize}
  \item[(i)] Arithmetical comprehension.
  \item[(ii)] $\forall {\fX}\;[\WO(\fX)\rightarrow \WO(2^{\fX})]$.
  \end{itemize}
  \end{thm}
  Another characterization from \cite{girard87}, Theorem 6.4.1, shows that arithmetical comprehension is equivalent to Gentzen's Hauptsatz (cut elimination)  for $\omega$-logic.
  Connecting statements of form (\ref{1}) to cut elimination theorems for
  infinitary logics will  be a major tool in this paper.
\subsection{$\mathbf{ACA}_0^+$  and $\varepsilon_{\fX}$}
 There are several more recent examples of such equivalences that have been proved by
 recursion-theoretic as well as proof-theoretic methods. The second example is a characterization in the same vein as (\ref{1}) for the theory $\mathbf{ACA}_0^+$
 in terms of the  $\xi\mapsto \varepsilon_{\xi}$ function.
  $\mathbf{ACA}_0^+$ stands for the theory $\mathbf{ACA}_0$ augmented
  by an axiom asserting that for any set $X$ the $\omega$ jump in $X$ exists. This theory is very  interesting as it is currently the weakest subsystem of $\Zz$ in which one has succeeded
to prove Hindman's Ramsey-type combinatorial theorem (asserting that finite colorings have infinite monochromatic sets closed under taking sums) and the Auslander/Ellis theorem of topological dynamics
(see \cite[X.3]{SOSA})
 The source of the pertaining representation system and its relativization is the Cantor normal form. 
 
 \begin{thm}  {\em (Cantor, 1897)}  For every ordinal 
$\beta>0$ there exist unique ordinals 
$\beta_0\geq\beta_1\geq\dots\geq\beta_n$ such that 
\begin{equation}\label{omega} \beta =
\omega^{\beta_0}+\ldots+\omega^{\beta_n}.\end{equation}
\end{thm}
The representation of $\beta$ in (\ref{epsilon}) is called
the
{\em Cantor normal form}.  We shall write
 $\beta\CNF\omega^{\beta_1}+\cdots+\omega^{\beta_n}$ to
convey that  $\beta_0\geq\beta_1\geq\dots\geq\beta_k$.
 
 $\varepsilon_0$ denotes the least ordinal 
$\alpha>0$ such that $
\beta<\alpha\,\;\Rightarrow\,\;\omega^{\beta}<\alpha$, or equivalently, 
the least ordinal  $\alpha$ such
that  $\omega^{\alpha}=\alpha$.
Every $\beta<\varepsilon_0$ has a Cantor normal form with
exponents $\beta_i<\beta$ and these exponents have Cantor
normal forms with yet again smaller exponents etc..  As this process
must terminate, ordinals $<\varepsilon_0$ can be effectively coded by
natural numbers.

 We state the result before introducing the functor ${\mathfrak X}\mapsto \varepsilon_{\fX}$.

\begin{thm}\label{MM}  {\em (Afshari, Rathjen \cite{AR}; Marcone, Montalb\'an \cite{MM})} Over ${\mathbf{RCA}}_0$ the following are equivalent:
 \begin{itemize}
  \item[(i)] $\mathbf{ACA}_0^+$.

  \item[(ii)] $\forall {\fX}\;[\WO(\fX)\rightarrow
\WO(\varepsilon_{\fX})]$.
  \end{itemize}
  \end{thm}
  
  Given a linear ordering $\fX=\langle X,<_X\rangle$, the idea for obtaining the new linear ordering $\varepsilon_{\fX}$ is to create a formal $\varepsilon$-number $\varepsilon_u$ for every $u\in X$ such that if $v<_X u$ then  $\varepsilon_v<_{\varepsilon_{\fX}}\varepsilon_u$, and in addition fill up the ``spaces" between these terms with formal Cantor normal forms.
   \paragraph{The ordering $<_{\varepsilon_{\fX}}$}
   \begin{deff}\label{epsilon}{\em 
  Let $\fX=\langle X,<_X\rangle$ be an ordering  where $X\subseteq{\mathbb N}$.
  $<_{\varepsilon_{\fX}}$ and its field $|\varepsilon_{\fX}|$ are inductively defined as follows:
\begin{enumerate}
\item $0\in |\varepsilon_{\fX}|$.
\item $\varepsilon_u\in |\varepsilon_{\fX}|$ for every $u\in X$, where
 $\varepsilon_u:=\langle 0,u\rangle$.
\item If $\alpha_1,\ldots,\alpha_n\in |\varepsilon_{\fX}|$, $n>1$ and
 $\alpha_n\leq_{\varepsilon_{\fX}}\ldots \leq_{\varepsilon_{\fX}}\alpha_1$,
  then $$\omega^{\alpha_1}+\ldots +\omega^{\alpha_n}\in |\varepsilon_{\fX}|$$ where
 $\omega^{\alpha_1}+\ldots +\omega^{\alpha_n}:=\langle 1,\langle\alpha_1,
\ldots,\alpha_n\rangle\rangle$.
\item  If $\alpha\in |\varepsilon_{\fX}|$ and $\alpha$ is not of the
form $\varepsilon_u$,
then $\omega^{\alpha}\in |\varepsilon_{\fX}|$,
where $\omega^{\alpha}:=\langle 2,\alpha\rangle$.

\item $0<_{\varepsilon_{\fX}}\varepsilon_u$ for all $u\in X$.
 \item $0<_{\varepsilon_{\fX}}\omega^{\alpha_1}+\ldots+ \omega^{\alpha_n}$ for all $\omega^{\alpha_1}+\ldots+ \omega^{\alpha_n}\in
 |\varepsilon_{\fX}|$.
 \item $\varepsilon_u<_{\varepsilon_{\fX}}\varepsilon_v$ if
 $u,v\in X$ and $u<_{_X}v$.
\item If $\omega^{\alpha_1}+\ldots +\omega^{\alpha_n}\in
 |\varepsilon_{\fX}|$, $u\in X$ and $\alpha_1<_{\varepsilon_{\fX}}\varepsilon_u$
 then  $\omega^{\alpha_1}+\ldots +\omega^{\alpha_n} <_{\varepsilon_{\fX}}
 \varepsilon_u$.
 \item If $\omega^{\alpha_1}+\ldots +\omega^{\alpha_n}\in
 |\varepsilon_{\fX}|$, $u\in X$, and $\varepsilon_u<_{\varepsilon_{\fX}}\alpha_1$ or $\varepsilon_u=\alpha_1$,
 then   $\varepsilon_u<_{\varepsilon_{\fX}} \omega^{\alpha_1}+\ldots +\omega^{\alpha_n}
 $.
 \item If $\omega^{\alpha_1}+\ldots+ \omega^{\alpha_n}$ and $\omega^{\beta_1}+\ldots+ \omega^{\beta_m}\in
 |\varepsilon_{\fX}|$, then
  \begin{eqnarray*}&&\omega^{\alpha_1}+\ldots+ \omega^{\alpha_n}
<_{\varepsilon_{\fX}} \omega^{\beta_1}+\ldots+ \omega^{\beta_m}\mbox{ iff} \\
 && n< m\,\wedge\,\forall i\leq n\;\alpha_i=\beta_i \;\;\mbox{ or }\\
 && \exists i\leq \min(n,m)[\alpha_i<_{\varepsilon_{\fX}}\beta_i\,\wedge\,
 \forall j<i\;\, \alpha_j=\beta_j]. \end{eqnarray*}

\end{enumerate}
 Let $\varepsilon_{\mathfrak X}=
  \langle  |\varepsilon_{\fX}|,<_{\varepsilon_{\fX}}\rangle$.}
  \end{deff}
One then proves (for instance in $\mathbf{RCA}_0$) that  $\varepsilon_{\mathfrak X}$ is linear ordering. If we denote by $\mathbb{LO}$ the category of linear orderings whose objects are all linear orderings and whose morphisms are the order-preserving maps between linear orderings, then ${\mathfrak X}\mapsto\varepsilon_{\mathfrak X}$ gives rise to an endofunctor of $\mathbb{LO}$:
For a morphism $f\co \mathfrak{X}\to \mathfrak{Y}$ define $\varepsilon(f)(t)$ for $t\in |\varepsilon_{\fX}|$ by replacing every expression $\varepsilon_u$ occurring in $t$ by
$\varepsilon_{f(u)}$. Then, $\varepsilon(f)\co \varepsilon_{\mathfrak X}\to \varepsilon_{\mathfrak Y}$. Moreover, the restriction of this functor to wellorderings is a dilator in the sense of Girard
\cite{gi} as it preserves pullbacks and direct limits (more about this  in subsection \ref{Gdilators}). 
  
   \paragraph{Ordinal representation
systems: 1904-1909}
At the beginning of the 20th century, mathematicians were intrigued by Cantor's continuum problem. Hardy in 1904 wanted
to ``construct'' a subset of $\mathbb{R}$ of size $\aleph_{1}$. In \cite{hardy}
  he gave explicit representations for all ordinals $<\omega^2$. Hardy's work then inspired O. Veblen who in his paper from 1908 \cite{veblen} found new normal forms for ordinals and succeeded in ``naming" 
 ordinals of an impressive chunk of the countable ordinals. In doing so, he also furnished proof theorists with the central idea for creating ordinal representation systems that were sufficient for much of their work until the middle 1960s. His $\varphi$-function was crucial in the work of S. Feferman and K. Sch\"utte who in the 1960s determined the limits of predicative mathematics with a notion of predicativity based on autonomous progressions of theories
(cf. \cite{f64,f68,sch64,sch65}). 
 
Veblen considered {\em continuous
increasing functions} on ordinals. Let $\ON$ be the class of ordinals. A (class) function
$f:\ON\to\ON$ is said to be {\em increasing} if $\alpha<\beta$
implies $f(\alpha)<f(\beta)$ and {\em continuous} (in the order
topology on $\ON$) if
$$f(\lim_{\xi<\lambda}\alpha_{\xi})=\lim_{\xi<\lambda}f(\alpha_{\xi})$$ holds for every
limit ordinal $\lambda$ and increasing sequence
$(\alpha_{\xi})_{\xi<\lambda}$.
$f$ is called {\em
normal} if it is increasing and continuous. By way of contrast,  the function $\beta\mapsto \omega+\beta$ is normal whereas
 $\beta\mapsto \beta+\omega$ is not since for instance 
 $\lim_{\xi<\omega}(\xi+\omega)=\omega$ but
  $(\lim_{\xi<\omega}\xi)+\omega=\omega+\omega$.

To these normal functions Veblen applied 
  two operations: \begin{itemize} \item[(i)]
{\em Derivation} \item[(ii)] {\em Transfinite Iteration}.
\end{itemize}
The {\em derivative} $f'$ of a function $f:\ON
\rightarrow \ON$ is the function which enumerates in increasing
order the solutions of the equation $$f( \alpha )= \alpha,$$ also
called the {\em fixed points} of $f$.
It is a fact of set theory that
 if $f$ is a normal function,
$$\{\alpha:\,f(\alpha)=\alpha\}$$ is a proper class and $f'$ will be
a normal function, too.

Using the two operations, Veblen defined a hierarchy of normal functions indexed along the ordinals.
 
\begin{deff}\label{veblenhierarchy}{\em 
 Given a normal function $f:\ON \rightarrow \ON$,
define a hierarchy of normal functions as follows:  \begin{enumerate} \item 
 $f_{0}\;=\;f $
 \item   $f_{\alpha +1}\;=\;{f_{\alpha}}'$
 \item  
$$f_{\lambda}(\xi)\; =\;\mbox{$\xi^{th}$ element of
}\bigcap_{\alpha<\lambda} \left\{\mbox{fixed points of } 
 f_{\alpha}\right\} \phantom{AA}\mbox{ for $\lambda$ limit}.$$ 
 \end{enumerate}
  Starting with the normal function $f(\xi)=\omega^{\xi}$,
the function $f_{\alpha}$ is usually denoted by 
$\varphi_{\alpha}.$
}\end{deff}
 The least ordinal $\gamma>0$ closed under  $\alpha,\beta\mapsto \varphi_{\alpha}(\beta)$, i.e.
the least ordinal $>0$ satisfying 
$(\forall\alpha,\beta<\gamma)\;\varphi_{\alpha}(\beta)<\gamma$
is the famous ordinal  $\Gamma_0$ that Feferman and Sch\"utte determined to be the least ordinal `unreachable'
by certain autonomous progressions of theories.

The two-place $\varphi$-function gives rise to a normal form theorem.

 \begin{thm}[$\vieh$ normal form]\label{4.27}
        For every additive principal\footnote{This means that $\alpha>0$ and $\forall \delta,\delta'<\alpha\; \delta+\delta'<\alpha$.} ordinal $\alpha$ there exist uniquely determined ordinals $\xi$ and $\eta$ such
        that $\alpha=\vieh_{\xi}(\eta)$ and $\eta<\alpha$.
        \end{thm}
        \prf See \cite[Theorem 13.12]{Sch77} or \cite[Theorem 5.27]{Ra12}. \qed

 The following comparison theorem encapsulates a procedure for recursively determining the order of $\varphi$-expressions, which can then be utilized to develop a representation system for the ordinals below $\Gamma_0$.

\begin{thm}[$\vieh$-comparison]\label{4.25}
         \begin{itemize}
         \item[(i)] $\vieh_{\alpha_1}(\beta_1)=\vieh_{\alpha_2}(\beta_2)$ holds iff
         one of the following conditions is satisfied:
         \begin{enumerate}
         \item $\alpha_1<\alpha_2$ and $\beta_1=\vieh_{\alpha_2}(\beta_2)$
        \item $\alpha_1=\alpha_2$ and $\beta_1=\beta_2$
        \item $\alpha_2<\alpha_1$ and $\vieh_{\alpha_1}(\beta_1)=\beta_2$.
        \end{enumerate}
       
         \item[(ii)] $\vieh_{\alpha_1}(\beta_1)<\vieh_{\alpha_2}(\beta_2)$ holds iff
         one of the following conditions is satisfied:
         \begin{enumerate}
         \item $\alpha_1<\alpha_2$ and $\beta_1<\vieh_{\alpha_2}(\beta_2)$
        \item $\alpha_1=\alpha_2$ and $\beta_1<\beta_2$
        \item $\alpha_2<\alpha_1$ and $\vieh_{\alpha_1}(\beta_1)<\beta_2$.
        \end{enumerate}
        \end{itemize}
      \end{thm}  
      \prf \cite[Theorems 13.9, 13.10]{Sch77} or \cite[Theorem 5.25]{Ra12}. \qed
      
      The representation system for $\Gamma_0$ can be relativized to any ordering $\mathfrak X$ by first introducing formal function symbols $\varphi_u$ for any $u\in X$ and secondly creating terms out of these
using also $0$ and formal addition $+$, and finally singling out the normal-forms. The crucial case in defining the ordering is the following:

\begin{eqnarray*} \varphi_u s <_{\varphi_{\mathfrak X}}\varphi_vt&\mbox{ iff }&
 u<_{_X}v\,\wedge\,s <_{\varphi_{\mathfrak X}} \varphi_v t \;\;\mbox{ or }\\
 && u=v\,\wedge\,s<_{\varphi_{\mathfrak X}}  t  \;\;\mbox{ or }\\ &&
v<_{_X}u\,\wedge\,\varphi_us <_{\varphi_{\mathfrak X}} t,
\end{eqnarray*}
   where all terms are assumed to be in normal form (for details see \cite{MM,rathjen-weiermann}). 
   
  $\varphi_{\mathfrak X}$ induces a functor on $\mathbb{LO}$ that characterizes $\mathbf{ATR}_0$.
\begin{thm}\label{zwei} (Friedman, unpublished; Marcone and Montalb\'an; Rathjen and Weiermann) Over ${\mathbf{RCA}}_0$ the following are equivalent:
 \begin{enumerate}
  \item  $\mathbf{ATR}_0$
  
  \item  $\forall {\fX}\;[\WO(\fX)\rightarrow
\WO(\varphi_{\fX})]$.
  \end{enumerate}
  \end{thm}
  Friedman's proof uses computability theory and also some proof theory.
  Among other things it uses a result which states that if
  $P\subseteq {\mathcal P}(\omega)\times {\mathcal P}(\omega)$ is arithmetic, then there is no
  sequence $\{A_n\mid n\in\omega\}$ such that
 (a) for every $n$, $A_{n+1}$ is the unique set such that
  $P(A_n,A_{n+1})$, and (b) for every $n$, $A'_{n+1}\leq_TA_n$.
  
  Of the two published proofs of the foregoing theorem, A. Marcone and A. Montalb\'an  employ tools from computability theory as their paper's title, {\em The Veblen function for computability
 theorists} \cite{MM}, clearly indicates.   A. Weiermann and the author of this paper use purely proof-theoretic means in \cite{rathjen-weiermann}.
 
 We would like to give the reader some insight into a proof strategy for showing theorems such as \ref{MM} and \ref{zwei}. However, we will be doing this by way of a different example. 
To state this result, it is convenient to introduce the notion of countable coded $\omega$-model.

\begin{deff}\label{omegamod}{\em Let $T$ be a theory in the language of second order arithmetic, $\mathcal{L}_2$. A {\em countable coded $\omega$-model of $T$} is a set $W\subseteq {\mathbb N}$, viewed as encoding the $\mathcal{L}_2$-model
  $${\mathbb M}=({\mathbb N},{\mathcal S},\in,+,\cdot,0,1,<)$$
  with ${\mathcal S}=\{(W)_n\mid n\in{\mathbb N}\}$ such that ${\mathbb M}\models T$ when the second order quantifiers are interpreted as ranging over ${\mathcal S}$ and the first order part is interpreted in the standard way
  (where $(W)_n=\{m\mid \paar{n,m}\in W\}$ with $\paar{\,,}$ being some primitive recursive coding function).

  If $T$ has only finitely many axioms, it is obvious how to express
  ${\mathbb M}\models T$ by just translating the second order quantifiers
  $QX\ldots X\ldots$ in the axioms by $Qx \ldots (W)_x\ldots$.
  If $T$ has infinitely many axioms,+ one needs to formalize Tarski's truth definition for ${\mathbb M}$.
This definition can be made in $\mathbf{RCA}_0$ as is shown in \cite{SOSA}, Definition II.8.3 and
  Definition VII.2. 

  We write $X\in W$ if $\exists n\;X=(W)_n$.
  }\end{deff}

The notion of countable coded $\omega$-model lends itself to alternative characterizations of Theorems \ref{MM} and \ref{zwei}.\footnote{The standard systems of reverse mathematics
 have the induction axiom $\forall X\,[0\in X\;\wedge\;\forall y\,(y\in X\to y+1\in X)\to \forall y\,y\in X]$ but aren't capable of deducing the induction scheme for all formulas of ${\mathcal L}_2$.
Note that  an $\omega$-model of such a theory will always satisfy the full induction scheme.}

\begin{thm}\label{alle} Over ${\mathbf{RCA}}_0$ the following are equivalent:
 \begin{itemize}
  \item[(i)] $\forall {\fX}\;[\WO(\fX)\rightarrow
\WO(\varepsilon_{\fX})]$ is equivalent to the statement that every set
 is contained in a countable coded $\omega$-model of $\mathbf{ACA}_0$.

  \item[(ii)] $\forall {\fX}\;[\WO(\fX)\rightarrow
\WO(\varphi{\fX}0)]$ is equivalent to the statement that every set
 is contained in a countable coded $\omega$-model of
$\Delta^1_1\mbox{-}{\mathbf{CA}}$ (or $\Sigma^1_1\mbox{-}{\mathbf{DC}}$).
  \end{itemize}
  \end{thm}
  \prf See \cite[Corollary 1.8]{H-Band}. \qed

  Whereas Theorem \ref{alle} has been  established independently by recursion-theoretic  and
  proof-theoretic methods, there is also a result that has a very involved proof and so far has only been shown by
  proof theory. It connects the $\Gamma$-function with the existence
  of countable coded $\omega$-models of ${\mathbf{ATR}}_0$.
  For this we  need to introduce the endofunctor $\mathfrak{X}\mapsto \Gamma_{\!\fX}$ of $\mathbb{LO}$. 
  The linear ordering $\Gamma_{\!\fX}$ is created out of $\mathfrak X$  by adding formal $\Gamma$-terms $\Gamma_u$ for every $u\in X$ with the stipulation that 
 $\Gamma_v<_{\Gamma_{\!\fX}}\Gamma_u$ if $v<_X u$,  and in addition one fills up the ``spaces" between these terms with formal sums and Veblen normal forms.
 The details can be found in \cite[Definition 2.5]{H-Band}.

   \begin{thm}\label{R} {\em (\cite[Theorem 1.4]{H-Band})}
    Over ${\mathbf{RCA}}_0$ the following are equivalent:
 \begin{itemize}
  \item[(i)]  $\forall {\fX}\;[\WO(\fX)\rightarrow
\WO(\Gamma_{\!\fX})]$.

  \item[(ii)] Every set
 is contained in a countable coded $\omega$-model of $\mathbf{ATR}_0$.
  \end{itemize}
  \end{thm}
  The tools from proof theory employed in the above theorems involve search trees and
  a plethora of cut elimination techniques for infinitary logic with ordinal bounds. 
The search tree techniques is a starting point that all proof-theoretic proofs of the theorems of this paper have in common. 
 It consists in producing the search or decomposition tree (in German ``Stammbaum") of a given formula. It proceeds by decomposing the formula according to its logical structure and amounts to applying
 logical rules backwards. This decomposition method has been employed by Sch\"utte \cite{sch56,sch60} to prove the completeness theorem for $\omega$-logic. It is closely related to the method of ``semantic tableaux"  of Beth \cite{beth} and the tableaux of
 Hintikka \cite{hintikka}. Ultimately, the whole idea derives from Gentzen \cite{gentzen}.

 \subsection{Proof idea  of (1)$\Rightarrow$(2) of Theorem \ref{R}}
 In it we shall use a simple result, namely that  $\mathbf{ATR}_0$ can be axiomatized
via a single sentence $\Pi^1_2$ sentence
$\forall X\,C(X)$
where $C(X)$ is $\Sigma^1_1$ (see \cite{SOSA}).

\begin{deff}{\em Let $\mathrm{L}_2$ be the language of second order arithmetic. We assume that $\mathrm{L}_2$ has relation symbols for primitive recursive relations. Formulas are generated from atomic  $R(t_1,\ldots,t_n)$ and negated atomic formulas $\neg R(t_1,\ldots,t_n)$ , where $R$ is a symbol an $n$-ary primitive recursive relation and $t_1,\ldots,t_n$ are numerical terms,\footnote{They are generated from numerical variables $x,y,z,\ldots$, the constants 0, 1, and the function symbols $+,\cdot$ for addition and multiplication, respectively.}
via the logical particles $\wedge,\vee, \exists x,\forall x,\exists X,\forall X$; so formulas are assumed to be in negation normal form.
 \begin{itemize}
 \item[(a)]  Let $U_0,U_1,U_2,\ldots$ be an enumeration of the free set variables of $\mathrm{L}_2$.
 \item[(b)]  For a closed term $t$, let $t^{^{\mathbb N}}$ be its numerical value. We shall assume that all predicate symbols of the language $\mathrm{L}_2$
are symbols for primitive recursive relations.  $\mathrm{L}_2$ contains
 predicate symbols for the primitive recursive relations of equality and inequality and possibly more (or all) primitive
 recursive relations. If $R$ is a predicate symbol in ${\mathrm L}_2$ we denote by $R^{^{\mathbb N}}$
the primitive recursive relation it stands for. If $t_1,\ldots,t_n$ are closed terms the formula $R(t_1,\ldots,t_n)$
($\neg R(t_1,\ldots,t_n)$)
is said to be {\em true} if $R^{^{\mathbb N}}(t_1^{^{\mathbb N}},\ldots,t_n^{^{\mathbb N}})$ is true (is false).
\item[(c)]  Henceforth a {\bf sequent} will be a finite list of ${\mathrm L}_2$-formulas {\it without free number variables}.
    Sequents will be denoted by upper case Greek letters. 
    
    Given sequents $\Gamma$ and $\Lambda$ and a formula $A$, we adopt the convention that $\Gamma,A,\Lambda$ denotes the sequent resulting from extending the list $\Gamma$ by $A$ and then extending it further by appending the list $\Lambda$. 
    \end{itemize}
\begin{itemize}
\item[(i)]  A sequent $\Gamma$ is {\bf axiomatic} if it satisfies at least  one of the following conditions:
\begin{enumerate}
 \item 
$\Gamma$ contains a
true {\bf literal}, i.e. a true formula of either form
$R(t_1,\ldots, t_n)$ or $\neg R(t_1,\ldots, t_n)$, where $R$ is a predicate symbol in $\mathrm{L}_2$ for a primitive recursive relation
and $t_1,\ldots,t_n$
are closed terms.
\item  $\Gamma$ contains formulas $s\in U$ and $t\notin U$
for some set variable $U$ and terms $s,t$ with $s^{^{\mathbb N}}=
t^{^{\mathbb N}}$.
 \end{enumerate}
 \item[(iv)]  A sequent is {\bf reducible} or a {\bf redex} if it is not axiomatic
 and contains a formula which is not a literal.
 \end{itemize}
}\end{deff}

 \begin{deff}{\em (Deduction chains in $\omega$-logic) 
 A {\bf deduction chain} is a finite string $$\Gamma_0;\Gamma_1;\ldots;\Gamma_k$$
of sequents $\Gamma_i$ constructed according to the following rules:
\begin{itemize}
\item[(i)]  $\Gamma_0\; = \;\emptyset$.
\item[(ii)]  $\Gamma_i$ is not axiomatic for $i<k$.
\item[(iii)]  If $i<k$ and $\Gamma_i$ is not reducible then
 $$\Gamma_{i+1}\;=\; \Gamma_i,\neg C(U_i).$$
   \item[(iv)]  Every reducible $\Gamma_i$ with $i<k$ is of the form
    $$\Gamma_i',E,\Gamma_i''$$
    where $E$ is not a literal and $\Gamma_i'$ contains only literals.\\
     $E$ is said to be the {\bf redex} of $\Gamma_i$.
\end{itemize}

Let $i<k$ and $\Gamma_i$ be reducible. $\Gamma_{i+1}$ is obtained from
 $\Gamma_i=\Gamma_i',E,\Gamma_i''$
     as follows:
     \begin{enumerate}
    \item  If $E\equiv E_0\vee E_1$ then
      $$\Gamma_{i+1}\;=\; \Gamma_i',E_0,E_1,\Gamma_i'',\neg C(U_i).$$
      \item   If $E\equiv E_0\wedge E_1 $ then
      $$\Gamma_{i+1}\;=\; \Gamma_i',E_j,\Gamma_i'',\neg C(U_i)$$
      where $j=0$ or $j=1$.
      \item   If $E\equiv \exists x\,F(x)$ then
      $$\Gamma_{i+1}\;=\; \Gamma_i',F(\bar{m}),\Gamma_i'',\neg C(U_i),E$$
      where $m$ is the first number such that $F(\bar{m})$ does not
occur in $\Gamma_0;\ldots;\Gamma_i$.

 \item If $E\equiv \forall x\,F(x)$ then
      $$\Gamma_{i+1}\;=\; \Gamma_i',F(\bar{m}),\Gamma_i'',\neg C(U_i)$$
      for some $m$.
      \item  If $E\equiv \exists X\,F(X)$ then
      $$\Gamma_{i+1}\;=\; \Gamma_i',F(U_{{m}}),\Gamma_i'',\neg C(U_i),E$$
      where $m$ is the first number such that $F(U_{{m}})$ does not
occur in $\Gamma_0;\ldots;\Gamma_i$.
 \item  If $E\equiv \forall X\,F(X)$ then
      $$\Gamma_{i+1}\;=\; \Gamma_i',F(U_{{m}}),\Gamma_i'',\neg C(U_i)$$
      where  $m$ is the first number such that $m\ne i+1$ and
$U_m$ does not
occur in $\Gamma_i$.

      \end{enumerate}
 }\end{deff}
  The set of deduction chains forms a tree ${\mathbb T}$ labeled with
strings of sequents. We will now consider two cases.
\paragraph{{\bf Case I:} ${\mathbb T}$ is not well-founded.}
   Then ${\mathbb T}$ contains an infinite path $\mathbb P$.
  Now define
 a set $M$ via
 \begin{eqnarray*} (M)_i &=& \{t^{^{\mathbb N}}\mid\mbox{$t\notin U_i$ occurs in ${\mathbb P}$}\}.\end{eqnarray*}
 Set ${\mathbb M}=({\mathbb N};\{(M)_i\mid i\in {\mathbb N}\},+,\cdot,0,1,<)$.

For a formula $F$,
 let $F\in {\mathbb P}$ mean that $F$ occurs in $\mathbb P$, i.e. $F\in \Gamma$ for some $\Gamma \in {\mathbb P}$.
 \\[1ex]
 {\bf Claim:} Under the assignment $U_i\mapsto (M)_i$ we have
 \begin{eqnarray*}\label{claim} F\in {\mathbb P} &\;\;\;\;\Rightarrow\;\;\;\;& {\mathbb M}\models
 \neg F.\end{eqnarray*}
 The Claim implies that ${\mathbb M}$ is an $\omega$-model of
$\mathbf{ATR}_0$. 
 

 \paragraph{{\bf Case II:} ${\mathbb T}$ is well-founded.}

We actually want to rule this  out. 
This is the place where the principle $$\forall {\fX}\;[\WO(\fX)\rightarrow
\WO(\Gamma_{\fX})]$$ in the guise of cut elimination for an
infinitary proof system  enters the stage.
Aiming at a contradiction, suppose that ${\mathbb T}$ is a well-founded tree.
Let $\fXO$ be the Kleene-Brouwer ordering on  ${\mathbb T}$.
Then $\fXO$ is a well-ordering.
In a nutshell, the idea is that a  well-founded ${\mathbb T}$ gives rise
 to a derivation of the empty sequent (contradiction) in the infinitary proof systems $\TUQ$
 from \cite{superjump}. This is were the really hard work lies and we have to stop here; details are in \cite{H-Band}.

  \section{Toward impredicative theories}
The proof-theoretic ordinal functions that figure in the foregoing theorems are all
  familiar from so-called predicative or meta-predicative  proof theory.
  Thus far a function from genuinely impredicative proof theory is missing.
  The first such function that comes to mind is of the Bachmann type \cite{bach}. We will shortly turn to it.
  
  Veblen in \cite{veblen} ventured very far into the transfinite in his 1908 paper, way beyond a representation system that incorporates the $\Gamma$-function.
  He extended the two-place $\varphi$-function first to an arbitrary finite
numbers of arguments, but then also to a transfinite
numbers of arguments, with the proviso that in, for example
$$\Phi_{f} ( \alpha_{0} , \alpha_{1} , \ldots , \alpha_{\eta} ),$$
only a finite number of the arguments $\alpha_{\nu}$ may be
non-zero.
 In this way Veblen singled out the ordinal {\blue $E(0)$}, which is often called the {\em big Veblen number}. $E(0)$ is
the least ordinal $\delta
> 0$ which cannot be named in terms of representations $$\Phi_{\ell} (
\alpha_{0} , \alpha_{1} , \ldots , \alpha_{\eta} )$$ with $\eta <
\delta$,
 and each $\alpha_{ \gamma } < \delta$.

\subsection{The Bachmann revelation.}
 In a paper published in 1950, Heinz Bachmann introduced a novel idea that allowed one to ``name" much larger ordinals than by the Veblen procedures. He had the amazing idea of using {\em uncountable ordinals} to keep track of diagonalizations over a hierarchies of functions. In more detail, he used the following steps:
 
\begin{itemize}
  \item 
  Define a set of
  ordinals {$\Bach$} closed under successor such that with each
  limit {$\lambda\in\Bach$} is associated an increasing sequence
{$\langle \lambda[\xi]:\,\xi<\tau_{\lambda}\rangle$} of
ordinals {$ \lambda[\xi]\in\Bach$} of length {\blue
$\tau_{\lambda}\in\Bach$} and
{$\lim_{\xi<\tau_{\lambda}}\lambda[\xi]=\lambda$}. \\[1ex]
 \item 
 Let $\OM$ be
  the {\em first uncountable ordinal}.
 A hierarchy of
functions $(\fiedi_{\alpha})_{\alpha\in\Bach}$ is then
obtained as follows:
\begin{eqnarray}\label{heinz}\label{bachmann} && \fiedi_0(\beta) \,=\, \omega^{\beta}\phantom{AAAA}
 \fiedi_{\alpha+1}\, =\, \left(\fiedi_{\alpha}\right)'\\ \nonumber
 && \fiedi_{\lambda}\;\mbox{ enumerates }\;
  \bigcap_{\xi<\tau_{\lambda}}(\mbox{Range of
  }\fiedi_{\lambda[\xi]})\phantom{aaa}\mbox{$\lambda$
  limit, $\tau_{\lambda}<\OM$} \\ \label{diag}
  && \fiedi_{\lambda}\;\mbox{ enumerates }\;
  \{\beta<\OM:\,\fiedi_{\lambda[\beta]}(0)=\beta\}
  \phantom{aaa}\mbox{$\lambda$
  limit, $\tau_{\lambda}=\OM$}.
  \end{eqnarray}
  \end{itemize}

 \paragraph{Distilling Bachmann's idea.} What makes Bachmann's approach rather difficult to deal with in proof theory is the requirement to  endow every limit ordinal with a fundamental sequence and 
referring to them when defining $\varphi_{\lambda}$, notably in the diagonalization procedure enshrined in (\ref{diag}). 
This layer of complication can be dispensed with though. What underpins the strength of Bachmann's approach can be described without the bookkeeping of fundamental sequences. 
We start by imagining a ``big" ordinal $\Omega$. What this means will become clearer as we go along, but definitely $\Omega$ should be an $\varepsilon$-number. 
 Using ordinals $<\Omega$ and $\Omega$ itself as building blocks, one then constructs further ordinals
using Cantor's normal form, i.e., if $\alpha_1\geq\ldots\geq \alpha_n$ have already been constructed, then
we build $\alpha:=\omega^{\alpha_1}+\cdots+\omega^{\alpha_n}$
provided that $\alpha>\alpha_1$. 
In this way we can build all ordinals $<\varepsilon_{\Omega+1}$, where the latter ordinal denotes the next 
$\varepsilon$-number after $\Omega$. 
Conversely, we can take any $\alpha< \varepsilon_{\Omega+1}$ apart, yielding smaller pieces as long as the exponents in its Cantor normal are smaller ordinals. This leads to the idea of support. More precisely define:  
\begin{deff}\label{support}{\em 
\begin{itemize} 
\item[(i)] $\EEO(0)=\emptyset$, $\EEO(\Omega)=\emptyset$.
\item[(ii)] $\EEO(\alpha)=\EEO(\alpha_1)\cup\cdots\cup\EEO(\alpha_n)$ if $\alpha=_{CNF}\omega^{\alpha_1}+\cdots+\omega^{\alpha_n}>
\alpha_1$.
\item[(iii)] {$\EEO(\alpha)=\{\alpha\}$} if $\alpha$ is an $\varepsilon$-number $<\Omega$.
\end{itemize} Note that $\EEO(\alpha)$ is a finite set.
 }\end{deff}
 To define something that is equivalent to what Bachmann achieves in (\ref{diag}),
the central idea is to devise an injective function {$$\varthetao\co \varepsilon_{\Omega+1} \to \Omega$$}
such that each $\varthetao(\alpha)$ is an $\varepsilon$-number.
 Think of $\varthetao$ as a {\em collapsing}, or better {\em projection} function in the sense of admissible set theory. 
For obvious reasons, $\varthetao$ cannot be order preserving, but the following can be realized: 
$$ \alpha<\beta\,\wedge\,\EEO(\alpha)<\varthetao(\beta) \;\to\; \varthetao(\alpha)<\varthetao(\beta).$$
 With the help of $\vartheta$-function\footnote{This function is a cousin of the $\theta$-function whose definition and properties are owed to Feferman, Aczel, Buchholz, and Bridge (see \cite[IX]{Sch77}) . $\vartheta$ was first defined in \cite{ra89} and used in \cite{rathjen-weiermann93}.} 
one obtains an ordinal representation system for the so-called {\em Bachmann ordinal}, also referred to as the {\em Bachmann-Howard ordinal}.\footnote{Bill Howard was looking for a description of the proof-theoretic ordinal of the theory of non-iterated inductive definitions. He was amazed and delighted to find it in  Bachmann's paper \cite{bach}.} 
\begin{deff}{\em We inductively define a set $\BHS$. \begin{itemize}
\item[(i)] $0\in \BHS$ and $\Omega\in\BHS$.
\item[(ii)] If $\alpha_1,\ldots,\alpha_n\in\BHS$, $\alpha_1\geq\cdots\geq\alpha_n$, then
$\omega^{\alpha_1}+\cdots+\omega^{\alpha_n}\in\BHS$. 
\item[(iii)] If $\alpha\in\BHS$ then so is $\varthetao(\alpha)$.
\end{itemize} 
$(\BHS,<)$ gives rise to an elementary ordinal representation system. Here $<$ denotes the restriction to $\BHS$. 
 
The {Bachmann-Howard ordinal} is the order-type of $\BHS\cap\Omega$.
}\end{deff}

 \subsection{Associating a Dilator with Bachmann.}
The Bachmann construction can be relativized to an arbitrary linear ordering $\mathfrak X$, as was shown in \cite{rathjen-pedro}, giving rise to a dilator $${\mathfrak X}\mapsto \vartheta_{\!_{\mathfrak X}}$$ and
the wellordering principle \begin{eqnarray}\label{WOPB}&&\WO({\mathfrak X})\Rightarrow \WO(\vartheta_{\!_{\mathfrak X}}).\end{eqnarray}
\begin{deff}\label{Bach}{\em \cite[2.6]{rathjen-pedro}
Again, let $\Omega$ be a ``big" ordinal.  Let $\mathfrak X$ be a well-ordering. With each $x\in X$ we associate a $\varepsilon$-number {$\Eb_x>\Omega$}. 
The set of ordinals $\BHSX$ is inductively defined as follows:
\begin{itemize}
\item[(i)] $0\in \BHSX$, $\Omega\in\BHSX$, and $\Eb_x\in\BHSX$ when $x\in X$.
\item[(ii)] If $\alpha_1,\ldots,\alpha_n\in\BHSX$, $\alpha_1\geq\cdots\geq\alpha_n$, then
$\omega^{\alpha_1}+\cdots+\omega^{\alpha_n}\in\BHSX$. 
\item[(iii)] If $\alpha\in\BHSX$ then so is $\varthetax(\alpha)$ and $\varthetax(\alpha)<\Omega$.
\end{itemize} 
To explain the ordering on $\BHSX$ one needs to extend $\EEO$: Let $\EEOX(\Eb_x)=\emptyset$.
 One then sets 
\begin{enumerate}
\item {$\Eb_x<\Eb_y\gdw x<_{\mathfrak X}y$}. 
\item $ \varthetax(\alpha)<\varthetax(\beta) \gdw ([ \alpha<\beta\,\wedge\,\EEOX(\alpha)<\varthetax(\beta)]\;\vee\;[\exists \gamma\in \EEOX(\beta)\;\varthetax(\alpha)\leq\gamma]).$
\end{enumerate} 
$(\BHSX,<)$ yields an  ordinal representation system elementary in $\mathfrak X$. 
}\end{deff}
The principle (\ref{WOPB}) turns ot to be equivalent over $\mathbf{RCA}_0$ to asserting that every set is contained in a countable coded $\omega$-model of the theory of {\em Bar Induction}, $\BI$,
whose formalization requires some preparations.

For a 2-place relation $\prec$ and an arbitrary formula $F(a)$ of $\mathcal{L}_2$ define
\begin{enumerate}
\item[]
$\text{Prog}(\prec,F):=(\forall x)[\forall y (y\prec x \rightarrow F(y))\rightarrow F(x)]$ (\emph{progressiveness})
\item[]$\text{\bf{TI}}(\prec,F):= \text{Prog}(\prec,F)\rightarrow\forall x F(x)$ (\emph{transfinite induction})
\item[] $\text{WF}(\prec):=\forall X\text{\bf{TI}}(\prec,X):=$ \newline
$\forall X(\forall x[\forall y (y\prec x \rightarrow y\in X))\rightarrow x\in X]\rightarrow \forall x[x\in X])$ (\emph{well-foundedness}).
\end{enumerate}
Let $\mathcal{F}$ be any collection of formulae of $\mathcal{L}_2$. For a 2-place relation $\prec$ we will write $\prec\in\mathcal{F}$, if $\prec$   is defined by a formula $Q(x,y)$ of $\mathcal{F}$ via $x\prec y:=Q(x,y)$.

\begin{dfn}{\em
$\mathrm{BI}$ denotes the bar induction scheme, i.e. all formulae  of the form
$$\text{WF}(\prec)\rightarrow\text{\bf{TI}}(\prec,F),$$
where $\prec$ is an arithmetical relation (set parameters allowed) and $F$ is an arbitrary formula of $\mathcal{L}_2$.

By $\BI$ we shall refer to the theory $\ACA_0+\mathrm{BI}$.
}\end{dfn}

The theorem alluded to above, due to P.F. Valencia Vizca\'ino and the author, conjectured in \cite[Conjecture 7.2]{rathjen-weiermann}, is the following:
 
 \begin{thm} {\em (Rathjen, Valencia Vizca\'ino \cite{rathjen-pedro})}
 Over ${\mathbf{RCA}}_0$ the following are equivalent:
 
 \begin{enumerate}
  \item  $\forall {\fX}\;[\WO(\fX)\rightarrow
\WO(\BHSX)]$.

\item  Every set is contained in a countable coded $\omega$-model of $\mathbf{BI}$.
  \end{enumerate}
   \end{thm} 

  \section{Towards a general theory of ordinal representations} We have seen several ordinal representation systems and their relativizations to an arbitrary well-ordering $\mathfrak X$.
Ordinal representation systems understood in this relativized way give rise to well-ordering principles: $$(*)\;\;\;\;\mbox{``if $\mathfrak X$ is a well-ordering then so is $\OR_{\mathfrak X}$"}$$ where
 $\OR_{\mathfrak X}$ arises from $\OR$ by letting $\mathfrak X$ play the role of the order-indiscernibles.
Any principle of the form $(*)$ is of syntactic complexity $\Pi^1_2$; thus cannot characterize stronger comprehensions such as
  $\Pi^1_1$-comprehension.\footnote{ $\Pi^1_1$-comprehension is syntactically of complexity $\Pi^1_3$. If it were equivalent to a $\Pi^1_2$ statement, say on the basis of $\mathbf{ATR}_0$, then it would follow that $\Pi^1_1\mbox{-}\mathbf{CA}_0$ proves its own consistency,  e.g. by \cite[Theorem VIII.5.12]{SOSA}.} We are therefore led to the idea of higher-order wellordering principles. To this end we require a general theory of ordinal representation systems.
 A crucial step towards a general theory via an axiomatic approach was taken by Feferman in  \cite{f68},  
 who singled out the property 
of {\em effective relative categoricity} as central to ordinal representation systems. 
Through the notion of relative categoricity he succeeded in crystallizing the algebraic aspect of ordinal representation systems by way of relativizing them to any set of order-indiscernibles, thereby in effect scaling them up to functors on the category of linear orders with order-preserving maps as morphisms. 
Below we recall his approach in \cite{f68}.
    
\subsection{Feferman's relative categoricity}  
      \begin{deff}{\em 
Let $f_1,\ldots, f_n$ be function symbols with $f_i$ having arity $m_i$.

The set of terms  $Tm(f_1,\ldots, f_n)$ is defined inductively as follows:
\begin{itemize}
\item[(i)] $\bar{0}\in Tm(f_1,\ldots, f_n)$;
 \item[(ii)] each variable $x_i$ is in $Tm(f_1,\ldots, f_n)$;
 \item[(iii)] if $t_1,\ldots,t_{m_i}\in Tm(f_1,\ldots, f_n)$, then
            so also       $$\mathsf{f}_i(t_1,\ldots,t_{m_i}).$$
            \end{itemize}
Let $\mathsf{Ord}$ denote the class of ordinals. Now suppose given a system of functions  $F_i\co \mathsf{Ord}^{m_i}\to \mathsf{Ord}$ (where $1\leq i\leq n$).  
By an {\em assignment} we mean any function $\sigma\co \omega\to    \mathsf{Ord}$. With each assignment $\sigma$ and any $t\in Tm(f_1,\ldots, f_n)$ is associated an ordinal $|t|_{\sigma}$, determined as follows.

\begin{enumerate}
\item $|\bar{0}|_{\sigma} = 0$;
\item $|x_k|=\sigma(k)$;
\item $|f_i(t_1,\ldots,t_{m_i})|_{\sigma}=F_i(|t_1|_{\sigma},\dots,|t_{m_i}|_{\sigma})$.
\end{enumerate} 
The system $\vec{F}$ is said to be {\bf replete} if for every ordinal $\alpha$ the closure of $\alpha$ under
$\vec{F}$ is an ordinal.

An ordinal $\kappa>0$ is said to be {\em $\vec{F}$-inaccessible} if whenever $\alpha_1,\ldots,\alpha_{m_i}<\rho$ then $F_i(\alpha_1,\ldots,\alpha_{m_i})<\rho$ holds for all $F_i$.
$\vec{F}$ is {\bf effectively relatively categorical} if, roughly speaking the order relation between any two terms $s(x_1,\ldots,x_k)$,   $t(x_1,\ldots,x_k)$ can be effectively determined from the ordering among $x_1,\ldots,x_k$
provided that these all represent $\vec{F}$-inaccessibles.
In particular  for assignments $\sigma,\tau$ into $\vec{F}$-inaccessibles it means that if
$$\forall i,j\in [1,\ldots,n](\sigma(x_i)<\sigma(x_j)\gdw \tau(x_i)<\tau(x_j))$$
then \begin{eqnarray*}|s(x_1,\ldots,x_k)|_{\sigma}<|t(x_1,\ldots,x_k)|_{\sigma} &\gdw& |s(x_1,\ldots,x_k)|_{\tau}<|t(x_1,\ldots,x_k)|_{\tau} \\
|t(x_1,\ldots,x_k)|_{\sigma}<|s(x_1,\ldots,x_k)|_{\sigma} &\gdw& |t(x_1,\ldots,x_k)|_{\tau}<|s(x_1,\ldots,x_k)|_{\tau}\\
|s(x_1,\ldots,x_k)|_{\sigma}=|t(x_1,\ldots,x_k)|_{\sigma} &\gdw& |s(x_1,\ldots,x_k)|_{\tau}=|t(x_1,\ldots,x_k)|_{\tau}.
\end{eqnarray*}
}\end{deff}
In his paper \cite{f70}, with the title {\em Hereditarily replete functionals over the ordinals},  Feferman  extended the foregoing notions also to finite type functionals over ordinals.

\subsection{Girard's dilators}\label{Gdilators}
A general theory of ordinal representation systems was also initiated by Girard who coined the notion of {\em dilator} \cite{gi}.

 Let  ${\Ord}$ be the category
whose objects are the ordinals and whose arrows are the strictly
increasing functions between ordinals. ${\Ord}$ enjoys two basic properties.

\begin{lem}\label{dilator1}\begin{itemize}
\item[(i)] $\Ord$ has pullbacks.
\item[(ii)] Every ordinal is a direct limit $\varinjlim(n_p,f_{pq})_{p\subseteq q\in I}$
 of a system of finite
 ordinals $n_p$ with $I$ being a set of finite sets of ordinals.
 \end{itemize}
 \end{lem}
 \prf  (i): Let $f\co \alpha \to \gamma$ and $g\co \beta \to \gamma$ be in $\Ord$. One easily checks that any $h\co \delta \to \gamma$
such that
\begin{eqnarray}\label{pb}\ran(h)&=&\ran(f)\cap \ran(g)\end{eqnarray}
is a pullback of $f\co \alpha \to \gamma$ and $g\co \beta \to \gamma$. Note that $h\co \delta \to \gamma$ is uniquely determined by (\ref{pb}) and that such an $h$ always exists:
Let  $\delta$ be the Mostowski collapse of $\ran(f)\cap \ran(g)$, $\clpse_{\ran(f)\cap \ran(g)}$ be the pertaining collapsing function, and $h$ be the inverse of the collapsing function. 

(ii): Fix an ordinal $\alpha$. Let $I$ be the collection of finite subsets of $\alpha$. $I$ is ordered by the inclusion relation, which is also directed,
i.e. for $s,t\in I$ there exists $q\in I$ such that $s\subseteq q$ and $t\subseteq q$. Let $s\in I$. Then $s=\{\alpha_0,\ldots,\alpha_{n-1}\}$ for uniquely determined  ordinals 
$\alpha_0<\ldots<\alpha_{n-1}$. Define $n_s:=n$ and $f_s\co n_s\to s$ by $f_s(i)=\alpha_i$. By design, $f_s$ is an order-preserving bijection between $n_s$ and $s$.
For $p\subseteq q\in I$ define $f_{pq}\co n_p\to n_q$ by $f_{pq}:= (f_q)^{-1}\circ f_p$. 
Thus, $f_{pp}$ is the identity on $n_p$, and if $p\subseteq q\subseteq r$, one has $f_{pr}= f_{qr}\circ f_{pq}$. 
As a result, $$\alpha=\varinjlim(n_p,f_{pq})_{p\subseteq q\in I}.$$
\qed


Of special interest are the endofunctors of this category that can be characterized by the two foregoing notions. 

\begin{deff}{\em  A {\em dilator} is an
endofunctor of the category  ${\Ord}$ preserving direct limits
 and pullbacks.}\end{deff}
 Dilators can be characterized in straightforward way.
 \begin{lem} For a  functor $F\co \Ord \to \Ord$, the following are equivalent: 
\begin{itemize}
\item[(i)] $F$ is a dilator.
\item[(ii)] $F$  has the following two properties:
 \begin{enumerate} \item  Whenever $f\co \alpha \to \gamma$, $g\co \beta \to \gamma$, $h\co \delta \to \gamma$, then
\begin{eqnarray}\label{pback}\ran(h)=\ran(f)\cap \ran(g) &\Rightarrow & \ran(F(h))=\ran(F(f))\cap \ran(F(g))\end{eqnarray}
\item For every ordinal $\alpha$ and $\eta<F(\alpha)$ there exists a finite ordinal $n$ and strictly increasing function $f\co n \to \alpha$
such that $\eta\in\ran(F(f))$.
\end{enumerate}
\end{itemize}
\end{lem}
\prf 
It follows from the proof of Lemma \ref{dilator1}(i) that preservation of pullbacks and condition (ii)(1) are equivalent. It remains to show that in $\Ord$ preservation of direct limits and condition 
(ii)(2) amount to the same. This can be gleaned from the proof of Lemma \ref{dilator1}(ii), however, as the details are not relevant for this paper, we just refer to \cite{gi} for the details. \qed

There is an easy but important consequence of the above, i.e. of Lemma \ref{dilator1}(ii).

\begin{lem} A dilator is completely determined by its behavior on the subcategory of finite ordinals and morphisms between them, $\Ord_{\omega}$. \end{lem}


At this point, it is perhaps not imeediately  gleanable  what the connection between ordinal representation systems and dilators might be. This is what we turn to next.

\paragraph{Denotation systems and dilators }
Girard (cf. \cite{gi,girard85,girardnorman85}) provided an alternative account of dilators that assimilates them more closely to ordinal representation systems and in particular to Feferman's notion
of relative categoricity.

\begin{deff}{\em Let $\ON$ be the class of ordinals and
 $F:\ON\to\ON$.
 A {\em denotation-system}
 for $F$ is a class $\mathcal D$ of ordinal {\em denotations}
 of the form $$(c;\alpha_0,\ldots,\alpha_{n-1};\alpha)$$
 together with an assignment $D:{\mathcal D}\to\ON$
 such that the following hold:
 
  \begin{enumerate} \item  If $(c;\alpha_0,\ldots,\alpha_{n-1};\alpha)$
 is in $\mathcal D$, then $\alpha_0<\ldots<\alpha_{n-1}<\alpha$
 and
  $D(c;\alpha_0,\ldots,\alpha_{n-1};\alpha)< F(\alpha)$.
 \item  Every $\beta<F(\alpha)$ has a unique denotation
  $(c;\alpha_0,\ldots,\alpha_{n-1};\alpha)$ in $\mathcal D$, i.e.\\
  $\beta=D(c;\alpha_0,\ldots,\alpha_{n-1};\alpha)$.
  \item  If $(c;\alpha_0,\ldots,\alpha_{n-1};\alpha)$ is a
  denotation and $\gamma_0<\ldots<\gamma_{n-1}<\gamma$, then
$(c;\gamma_0,\ldots,\gamma_{n-1};\gamma)$ is a denotation.
 \item  If $D(c;\alpha_0,\ldots,\alpha_{n-1};\alpha)\leq
 D(d;\alpha_0',\ldots,\alpha_{m-1}';\alpha)$,
 $\gamma_0<\ldots<\gamma_{n-1}<\gamma$,
  $\gamma_0'<\ldots<\gamma_{m-1}'<\gamma$, and for all
  $i<n$ and $j<m$, $$\alpha_i\leq\alpha_j'\Leftrightarrow
  \gamma_i\leq \gamma_j' \;\;\mbox{ and }\;\; \alpha_i\geq\alpha_j'\Leftrightarrow
  \gamma_i\geq \gamma_j'$$ then
  $$D(c;\gamma_0,\ldots,\gamma_{n-1};\gamma)\leq
 D(d;\gamma_0',\ldots,\gamma_{m-1}';\gamma).$$
\end{enumerate}

In a denotation $(c;\alpha_0,\ldots,\alpha_{n-1};\alpha)$,
    $c$ is called the {\em index}, $\alpha$ is the {\em parameter} and
     $\alpha_0,\ldots,\alpha_{n-1}$ are the {\blue coefficients} of
     the denotation. If $\beta=D(c;\alpha_0,\ldots,\alpha_{n-1};\alpha)$,
     the index $c$ represents some `algebraic' way of
     describing $\beta$ in terms of the ordinals
     $\alpha_0,\ldots,\alpha_{n-1},\alpha$.
     }\end{deff}

 \begin{lem}  Every
denotation system $D$ induces a dilator $F_{\!_D}$ by letting
 $F_{\!_D}(\alpha)$ be the least ordinal $\eta$ that does not have
  a denotation of the form
  $D(c;\alpha_0,\ldots,\alpha_{n-1};\alpha)$, and for any arrow
   $f:\alpha\to\delta$ of the category $\Ord$ letting
   $$F_{\!_D}(f):F_{\!_D}(\alpha)\to F_{\!_D}(\delta)$$ be defined by
    $$F_{\!_D}(f)(D(c;\alpha_0,\ldots,\alpha_{n-1};\alpha))\;:=\;
D(c;f(\alpha_0),\ldots,f(\alpha_{n-1});\delta).$$
\end{lem}
The converse is also true. 
 
\begin{lem} To every dilator $F$ 
  one can assign a 
 denotation system $D_{\!_F}$ such that $\beta<F(\alpha)$ is denoted by
$$(\gamma;\alpha_0,\ldots,\alpha_{n-1};\alpha)$$
where $n$ is the least finite ordinal such that there exists a morphism
$f\co n\to \alpha$ with $\beta\in \ran(F(f))$. Moreover, $\gamma<F(n)$ is uniquely determined by $F(f)(\gamma)=\beta$, and
$\alpha_0=f(0),\ldots,\alpha_{n-1}=f(n-1)$.
\end{lem}
\prf $F$ being a dilator, we know that such an $f$ exists. $f$ is also uniquely determined since $n$ is chosen minimal: for if $g\co n\to \alpha$ also satisfied 
$\beta\in \ran(F(g))$, then their joint pullback $h\co m\to \alpha$ would satisfy $\beta\in \ran(h)$ as well, so that $n=m$ and $\ran(h)=\ran(f)=\ran(g)$, meaning $f=g$.
\qed

 \section{Higher order wellordering principles}
 We already noted that a statement of the form $\WOP(f)$ is $\Pi^1_2$ and therefore cannot be equivalent
 to  a theory whose axioms have an essentially higher complexity, like for instance $\Pi^1_1$-comprehension.
 The idea thus is that to get equivalences one has to climb up in the type structure. Given a functor
$$F:(\LOO\to\LOO)\to(\LOO\to\LOO),$$ where $\LOO$ is the class of linear orderings, we consider the statement:\footnote{$\WOPP$ being an acronym for the German ``{\bf H}\"oherstufiges {\bf W}ohl{\bf o}rdnungs{\bf p}rinzip''.}
$$\WOPP(F):\phantom{AAA}\forall f\in (\LOO\to\LOO)\;[\WOP(f)\to\WOP(F(f))].$$
There is also a variant of $\WOPP(F)$ which should basically encapsulate the same ``power''.
Given a functor
$$G:(\LOO\to\LOO)\to\LOO$$ consider the statement:
$$\WOPPP(G):\phantom{AAA}\forall f\in (\LOO\to\LOO)\;[\WOP(f)\to\WO(G(f))].$$

 It was conjectured in \cite{montalban} that a principle of the form $\WOPP(F)$ might be equivalent to $\Pi^1_1$-comprehension. In \cite{rathjen-chicago} is was claimed that for a specific functor
${\mathcal B}:(\LOO\to\LOO)\to\LOO$, $\WOPPP({\mathcal B})$ is equivalent to $\Pi^1_1$-comprehension. ${\mathcal B}$ is a functor that takes an arbitrary dilator $F$ as input and returns an ordinal representation system; this 
 amounts to combining the Bachmann procedure with the closure under $F$.  \cite{rathjen-chicago} also adumbrated the steps and a proof strategy for this result.

  \subsection{Bachmann meets a dilator}
The idea to combine the Bachmann construction with an arbitrary dilator one finds in \cite{montalban} and \cite{rathjen-chicago}. The details were worked out in
 Anton Freund's thesis \cite{freund-thesis}.\footnote{Of course, Definition \ref{Bach}, which harks back to \cite{rathjen-pedro}, set an important precedent.}
  
  
\begin{deff}\label{Bachdil}{\em Again, let $\Omega$ be a ``big" ordinal.  Let $\mathcal D$ be a system of denotations.
 
\begin{itemize}
\item[(i)] $0\in \BHSD$, $\Omega\in\BHSD$.
\item[(ii)] If $\alpha_1,\ldots,\alpha_n\in\BHSD$, $\alpha_1\geq\cdots\geq\alpha_n$, then
$\omega^{\alpha_1}+\cdots+\omega^{\alpha_n}\in\BHSD$.
\item[(iii)] If $\alpha\in\BHSD$ then so is $\varthetad(\alpha)$ and $\varthetad(\alpha)<\Omega$.

\item[(iv)] If  $\alpha_0,\ldots,\alpha_{n-1}\in\BHSD$, $\alpha_0<\cdots<\alpha_{n-1}<\Omega$, and 
$\sigma=(c;0,\ldots,n-1;n)\in \mathcal{D}$, then
$$\EDD^{\sigma}_{\alpha_0,\ldots,\alpha_{n-1}}\in \BHSD\mbox{ and }\Omega<\EDD^{\sigma}_{\alpha_0,\ldots,\alpha_{n-1}}.$$ 
\end{itemize} 
 To explain the ordering on $\BHSD$ one needs to extend $\EEO$:  Let  $$\EEOD(\EDD^{\sigma}_{\alpha_0,\ldots,\alpha_{n-1}})=\{\alpha_0,\ldots,\alpha_{n-1}\}.$$
The ordering between $\varthetad$-terms is defined as before:
$$ \varthetad(\alpha)<\varthetad(\beta) \gdw ([\alpha<\beta\,\wedge\,\EEOD(\alpha)<\varthetad(\beta)]\;\vee\;[\exists \gamma\in \EEOD(\beta)\, \varthetad(\alpha)\leq \gamma]).$$ 
How do we compare terms of the form  $\EDD^{\sigma}_{\alpha_0,\ldots,\alpha_{n-1}}$ and 
 $\EDD^{\tau}_{\beta_0,\ldots,\beta_{m-1}}$?
To this end  let $\sigma=(c;0,\ldots,n-1;n)$ and $\tau=(e;0,\ldots,m-1;m)$; also let $k$ be the number of elements of  $\{\alpha_0,\ldots,\alpha_{n-1},\beta_0,\ldots,\beta_{m-1}\}$.
 In a first step, define $f\co n\to k$ and $g\co m\to k$
such that
\begin{eqnarray*} f(i)< g(j) &\gdw & \alpha_i<\beta_j \\
f(i)=g(j) &\gdw & \alpha_i=\beta_j.\end{eqnarray*} 
Then $$\EDD^{\sigma}_{\alpha_0,\ldots,\alpha_{n-1}} \;<\;\EDD^{\tau}_{\beta_0,\ldots,\beta_{m-1}}$$ if and only if
 $$D(c;f(0),\ldots,f(n-1);k) \;<\;D(e;g(0),\ldots,g(m-1);k).$$
 }\end{deff}
It is rather rather instructive to see why $\BHSD$ is a wellordering.We shall prove that with the help of strong comprehension. 

\begin{thm}\label{WOBD} {\em ($\Pi^1_1\mbox{-}\mathbf{CA}_0$)}
 For every dilator $\mathcal D$, $(\BHSD,<)$ is a wellordering.
\end{thm}
\prf Let $I$ be the well-founded part of $\BHSD\cap\Omega:=\{\alpha\in\BHSD\mid \alpha<\Omega\}$. $\Pi^1_1$-comprehension ensures that $I$ is a set. Then ${\mathfrak I}:=(I,<\restriction I\times I)$ is a wellordering. 
Now let $${M}\,:=\,\{\beta\in\BHSD\mid \EEOD(\beta)\subseteq I\}\;\mbox{ and }\;\mathfrak{M}\,:=\, (M,<\restriction M\times M).$$
We claim that $\mathfrak{M}$ is a wellordering, too. To see this, note that the set of terms of the form $\EDD^{\sigma}_{\alpha_0,\ldots,\alpha_{n-1}}$ with $\alpha_0,\ldots,\alpha_{n-1}\in I$ is wellordered as they are in one-one and order preserving correspondence with the denotations $D(c;\alpha_0,\ldots,\alpha_{n-1};\mathfrak{I})$ where $\sigma=(c;0,\ldots,n-1;n)$. Now, $M$ is obtained from $I$ by adding the latter terms as well as $\Omega$, and then  closing off under
$+$ and $\xi\mapsto \omega^{\xi}$ (or more precisely, Cantor normal forms). It is known from Gentzen's \cite{gentzen43} that this last step preserves wellorderedness.\footnote{For details, consult  \cite[VIII.2]{Sch77},  \cite[9.6.2]{pohlers} or \cite[Section 4]{rathjen2005}.}

Thus, if we can show that  $\BHSD$ is the same as $\mathfrak M$ we are done. For this it is enough to establish closure of $M$ under $\varthetad$.
$$\mathsf{Claim:}\;\;\; \forall \alpha \in M\;\varthetad(\alpha)\in M.$$
The Claim is proved by induction on $\alpha$. So assume $\alpha\in M$ and  $$(*)\;\;\;\forall \beta\in M\;[\beta<\alpha\to \varthetad(\beta)\in I].$$
To be able to conclude that $\varthetad(\alpha)\in I$, it suffices to show that $\delta\in I$ holds for all $\delta< \varthetad(\beta)$. 
We establish this via a subsidiary induction on the syntactic complexity of $\delta$ (so this is in effect an induction on naturals). If $\delta=0$ this is immediate, and if $\delta=\omega^{\delta_1}+\ldots+\omega^{\delta_n}$ with
$\delta>\delta_1\geq \ldots\geq\delta_n$ this follows from the subsidiary induction hypothesis and the fact that $I$ is closed under $+$ and $\xi\mapsto \omega^{\xi}$.

This leaves the case that $\delta=\varthetad(\rho)$ for some $\rho$. Since $\varthetad(\rho)< \varthetad(\alpha)$ there are two subcases to consider.

{\em Case 1}: $\rho<\alpha$ and  $\EEOD(\rho)<\varthetad(\alpha)$. By the subsidiary induction hypothesis, $\EEOD(\rho)\subseteq I$ as the syntactic complexity of terms in $\EEOD(\rho)$ is not bigger than that of $\rho$ and thus smaller than that of $\delta$. Therefore, $\rho\in M$, and since $\rho<\alpha$ it follows from $(*)$ that $ \delta=\varthetad(\rho)\in I$.

{\em Case 2}: There exists $\xi\in \EEOD(\alpha)$ such that $\delta\leq \xi$. $\alpha\in M$ entails that $\xi\in I$, and thus $\delta\in I$.
\qed

  \begin{thm}\label{haupt} {\em (Freund)} The following are equivalent over $\mathbf{RCA}_0$:
  \begin{itemize}
  \item[(i)] $\Pi^1_1$-comprehension.
  \item[(ii)] For every denotation system $\mathcal D$, $\BHSD$ is a wellordering.
  \end{itemize}
  \end{thm}
  A detailed proof was given in Anton Freund's thesis \cite{freund-thesis} and in his papers \cite{freund0,freund1,freund2}.
  Below, however, we will be presenting another proof, partly going back to \cite{rathjen-chicago}, since the purposes here is to convey the intuitions behind the result to the reader and present them without employing 
  category-theoretic terminology. 
  
The direction (i)$\Rightarrow$(ii) of Theorem \ref{haupt} is taken care of by Theorem \ref {WOBD}. 
 Next,  we will be concerned with giving a proof sketch for the direction $(*)$  (ii)$\Rightarrow$(i). 
 For $(*)$  it is enough to show that every set $\bar{Q}\subseteq \mathbb{N}$ is contained in a countable coded model ${\mathfrak M}=(M,E)$ of Kripke-Platek set theory (with Infinity), where the interpretation $E$ of the elementhood relation is well-founded. This is so because a  $\Pi^1_1$-definable class of naturals, $\{n\in\mathbb{N}\mid F(n,\bar{Q})\}$, with second order parameter $\bar{Q}$ is $\Sigma_1$-definable over $\mathfrak M$ (see \cite[IV.3.1]{ba}), and therefore a set in the background theory. Important techniques include Sch\"utte's method of search trees \cite{sch56} 
 for $\omega$-logic but adjusted to $\alpha$-logic for any ordinal (or rather well-ordering) $\alpha$ and proof calculi similar to the one used for the ordinal analysis of Kripke-Platek set theory (originally due to J\"ager \cite{j80,j82}). It is not necessary, though, to use the machinery for the ordinal analysis of $\mathbf{KP}$. Instead, one can also directly 
 employ the older techniques from ordinal analyses of the formal theory  of non-iterated arithmetic inductive definitions as presented in \cite{pohlers81} and \cite{pohlers}.
 Technically, though, the easiest approach appears to be  to work with an extension of Peano arithmetic  with inductive definitions and an extra sort of ordinals to
express the stages of an inductive definition as in J\"ager's $\mathbf{PA}_{\Omega}$ from  \cite{j93}. 

\begin{deff}{\em In what follows, we fix a set $\bar{Q}\subseteq\mathbb{N}$ and a well-ordering $\mathfrak{X}=(X,<_X)$ with $X\subseteq \mathbb{N}$. 
Let $L_1$  be the first-order language of arithmetic with number
variables $x, y, z, \ldots$  (possibly with subscripts), the constant $0$, as well as function
and relation symbols for all primitive recursive functions and relations, and a unary predicate symbol $Q$.  Let $ L_1(P)$
be the extension of $ L_1$  by a further ``fresh"  unary relation symbol $P$. 
The {\em atomic formulas} of  $ L_1(P)$ are of the form $Q(t)$,  $P(t)$ and $ R(s_1,\ldots,s_k)$, where $t,s_1,\ldots,s_k$ are terms and $R$ is a relation symbol for a primitive recursive $k$-ary relation. The {\em literals} of $ L_1(P)$ are the atomic formulas and their negations  $\neg Q(t)$, $\neg P(t)$ and $ \neg R(s_1,\ldots,s_k)$. The {\em formulas} of  $ L_1(P)$ are then generated from
the literals via the logical connectives $\wedge,\vee$ and the quantifiers $\forall x$ and $\exists x$. Note that the formulas are in negation normal form. As per usual, negation is a defined operation, using deletion of double negations and de Morgan's laws.
A formula is said to be {\em $P$-positive} if it contains no occurrences of the form $\neg P(t)$.  A $P$-positive formula in which at most the variable $x$ occurs free will be called an {\em inductive operator form}; expressions $A(P,x)$ are meant to range over such forms.

$L_1$ will be further extended  to a new  language $\LQ$ by adding a new sort of
ordinal variables $\alpha,\beta,\gamma,\ldots$ (possibly with subscripts), a new binary relation
symbol $<$ for the less relation on the ordinals\footnote{The usual ordering on the naturals will be notated by $<_N$.} and a binary relation symbol $P_A$ for
each inductive operator form $A(P, x)$.\footnote{Observe that the place holder predicate $P$ of  $ L_1(P)$ is not part of the language $\LQ$.}  For each element $u$ of the wellordering $\mathfrak{X}$ we also introduce an ordinal constant $\alpha_u$ into the language $\LQ$.
Ordinal constants inherit an ordering from $\mathfrak X$, namely $\alpha_u$ is considered to be smaller than $\alpha_v$ if $u<_Xv$. In future,  when we talk about the least ordinal constant satisfying a certain condition we are referring to that ordering. 

Since subscripts are often a nuisance to humans, we will use overlined ordinal variables, i.e.  $\bar{\alpha},\bar{\beta},\ldots$,  as metavariables for these constants.

The number terms $ s, t, \ldots$   of $\LQ$ are the number
terms of $ L_1$; the ordinal terms of $\LQ$ are the ordinal variables and constants. The formulas of $\LQ$ 
 are inductively generated as follows:
\begin{enumerate}
\item  If $R$ is a $k$-ary relation symbol of $L_1$ and $s_1,\ldots,s_k$ are number terms, then $R(s_1,\ldots,s_k)$  is an atomic formula of $\LQ$.

\item $\mathfrak{a}<\mathfrak{b}$, $\mathfrak{a}=\mathfrak{b}$  and $P_A(s,\mathfrak{a})$ 
are  atomic formulas of $\LQ$ whenever $\mathfrak{a}$ and 
$\mathfrak{b}$ are ordinal terms and $s$ is a number term.

\item Negations of atomic formulas of $\LQ$ are formulas of $\LQ$. The latter together with the atomic formulas make up the {\em literals} of $\LQ$.

\item If $B,C$ are formulas of $\LQ$, then so are $B\wedge C$ and $B\vee C$.
\item If $F(x)$ is a formula of $\LQ$, then so are $\exists xF(x)$ and $\forall xF(x)$. 
\item If $G(\alpha)$ is a formula of $\LQ$, then so are $\exists\alpha G(\alpha)$ and $\forall \alpha F(\alpha)$.
\item  If $G(\alpha)$ is a formula of $\LQ$ and $\mathfrak{b}$ is an ordinal term, then so are $(\exists\alpha<\mathfrak{b})\, G(\alpha)$ and $(\forall \alpha<\mathfrak{b})\, F(\alpha)$.
\end{enumerate}
As per usual, $B\to C$ is an abbreviation for $\neg B\,\vee\,C$  (with $\neg$ being the obvious defined operation of negation), and $B\leftrightarrow C$ abbreviates $(B\to C)\,\wedge\,(C\to B)$.

In future,  we write $P^{\mathfrak{a}}_A(s)$ for $P_A(s,\mathfrak{a})$. We also write $P^{<\mathfrak{a}}_A(s)$ for $(\exists\beta<\mathfrak{a})\,P^{\beta}_A(s)$ and $P_A(s)$ for $\exists\alpha\,P^{\alpha}_A(s)$.

\comment{We also distinguish two important classes of formulas. The {\em $\DeltaO$-formulas} are obtained from the atomic and negated atomic formulas by closing off under $\wedge$, $\vee$, number quantifiers, and  bounded ordinal quantifiers  $\exists\alpha<\mathfrak{b}$ and $\forall \alpha<\mathfrak{b}$ with $\mathfrak{b}$ an ordinal term.

The {\em $\SigmaO$-formulas} are those obtained from the atomic and negated atomic formulas by closing off under $\wedge$, $\vee$, number quantifiers,  bounded ordinal quantifiers,
and unbounded existential ordinal quantifiers $\exists \alpha$.}
}\end{deff}

\begin{deff}\label{Axiome}{\em The axioms of $\PAO$ fall into several groups:
\begin{enumerate}
\item True literals: These are of three forms. 
\begin{enumerate}
\item Let $R$ be a symbol for a $k$-ary primitive recursive relation and $s_1,\ldots,s_k$ are closed number terms.   $R(s_1,\ldots,s_k)$ is a true literal if $\mathcal{R}(s^{\mathbb N}_1,\ldots, s^{\mathbb N}_1)$ is true in the naturals, where $\mathcal{R}$ denotes the predicate pertaining to $R$ and  $s^{\mathbb{N}}_i$ denotes the value
of the closed number term $s_i$ in the standard interpretation; if $R(s_1,\ldots,s_k)$ is false under the standard interpretation, then  $\neg R(s_1,\ldots,s_k)$  is a true literal.

\item $Q(s)$ is a true literal if  $s^{\mathbb{N}}_i\in \bar{Q}$, and $\neg Q(t)$ is a true literal if  $t^{\mathbb{N}}_i\notin \bar{Q}$, where $s,t$ are closed numeral terms.

\item Let $\alpha_u,\alpha_v$ are ordinal constants with $u,v\in \mathfrak{X}$.  $\alpha_u<\alpha_v$ is a true literal if $u<_Xv$, and otherwise $\neg\,\alpha_u<\alpha_v$ is a true literal.
$\alpha_u=\alpha_u$ is a true literal whereas $\neg \alpha_u=\alpha_v$ is a true literal if $u\ne v$.
\end{enumerate}
\item Stage axioms: These are all sentences of the form $\forall\alpha\,\forall x\,[A(P_A^{<\alpha},x)\to P_A^{\alpha}(x)]$.
\item Reflection: These are all sentences of the form $\forall x\,[A(P_A,x)\to \exists\alpha\,P^{\alpha}_A(s)]$.
\item Let $A_0,A_1,\ldots$ be a fixed enumeration of all stage and reflection axioms. Note that all of these are closed formulas (sentences) and that none of them contains any
ordinal constants.
\end{enumerate}

}\end{deff}

   \subsection{Deduction chains in $\PAO$}\label{DedC}

Sequents are finite lists of sentences of $\LQ$; they will be notated by upper case Greek letters.\footnote{As before, $\Gamma,A$ denotes the sequent obtained from $\Gamma$ via appending the formula $A$ to the list $\Gamma$. Similarly, $\Gamma,A_1,\dots,A_r,\Lambda$
stands for the list $C_1,\ldots,C_n,A_1,\dots,A_r,D_1,\ldots,D_q$ if $\Gamma$ is the list  $C_1,\ldots,C_n$ and $\Lambda$ is the list  $D_1,\ldots,D_q$.}

\begin{deff}{\em \begin{itemize}
\item[(i)]  A sequent $\Gamma$ is {\em  axiomatic} if it  contains a
true {\em literal}.\footnote{Note that if $\Gamma$ contains formulas $R(\vec{s}\,)$ and $\neg R(\vec{t}\,)$, where $s_i^{\mathbb{N}}=t_i^{\mathbb{N}}$ and $R$ is a symbol for a primitive recursive 
relation or equals ${Q}$, then $\Gamma$ is axiomatic.}
 \item[(ii)]  A sequent is {\em  reducible} or a {\em redex} if it is not axiomatic
 and contains a formula which is not a literal.
 \end{itemize}
}
\end{deff}

\begin{deff}{\em 
 A {\em $\PAO$-deduction chain} is a finite string $$\Gamma_0;\Gamma_1;\ldots;\Gamma_k$$
of sequents $\Gamma_i$ constructed according to the following rules:
\begin{itemize}
\item[(i)]  $\Gamma_0\; = \;\emptyset$.
\item[(ii)]  $\Gamma_i$ is not axiomatic for $i<k$.
\item[(iii)]  If $i<k$ and $\Gamma_i$ is not reducible then
 $$\Gamma_{i+1}\;=\; \Gamma_i,\neg A_i$$
where $A_i$ is the $i$-th formula in the list from Definition \ref{Axiome}(4).
   \item[(iv)]  Every reducible $\Gamma_i$ with $i<k$ is of the form
    $$\Gamma_i',E,\Gamma_i''$$
    where $E$ is not a literal and $\Gamma_i'$ contains only literals.
     $E$ is said to be the {\em redex} of $\Gamma_i$.

 Let $i<k$ and $\Gamma_i$ be reducible. $\Gamma_{i+1}$ is obtained from
 $\Gamma_i=\Gamma_i',E,\Gamma_i''$
     as follows:
     \begin{enumerate}
    \item  If $E\equiv E_0\vee E_1$ then
      $$\Gamma_{i+1}\;=\; \Gamma_i',E_0,E_1,\Gamma_i'',\neg A_i.$$
      \item   If $E\equiv E_0\wedge E_1 $ then
      $$\Gamma_{i+1}\;=\; \Gamma_i',E_j,\Gamma_i'',\neg A_i$$
      where $j=0$ or $j=1$.
      \item   If $E\equiv \exists x\in {\mathbb N}\,F(x)$ then
      $$\Gamma_{i+1}\;=\; \Gamma_i',F(\bar{m}),\Gamma_i'',\neg A_i,E$$
      where $m$ is the first number such that $F(\bar{m})$ does not
occur in $\Gamma_0;\ldots;\Gamma_i$.\footnote{There is the slightly irksome possibility that $x$ does not occur free in $F(x)$. Then $m$ can be any number for it doesn't matter
which number one choses.}
 \item  If $E\equiv \forall x\in{\mathbb N}\,F(x)$ then
      $$\Gamma_{i+1}\;=\; \Gamma_i',F(\bar{m}),\Gamma_i'',\neg A_i$$
      for some $m$.
      \item If $E\equiv \exists \xi\,F(\xi)$ then
      $$\Gamma_{i+1}\;=\; \Gamma_i',F(\bar{\delta}),\Gamma_i'',\neg A_i,E$$
      where the ordinal constant $\bar{\delta}$ is picked as follows. If there occurs an ordinal constant $\bar{\gamma}$ in  $\Gamma_0;\ldots;\Gamma_i$ such that $F(\bar{\gamma})$ does not
occur in $\Gamma_0;\ldots;\Gamma_i$ then let $\bar{\delta}$ be the least such (with leastness understood in the sense of $\mathfrak X$).
If, however, for all ordinal constants  $\bar{\gamma}$  in  $\Gamma_0;\ldots;\Gamma_i$ the formula $F(\bar{\gamma})$ already
occurs in $\Gamma_0;\ldots;\Gamma_i$, then $\bar{\delta}$ can be any ordinal constant.\footnote{One might be tempted to decree that $\bar{\delta}$ be the first ordinal constant (in the sense of $\mathfrak X$)
that does not occur in $\Gamma_0;\ldots;\Gamma_i$. However, this is problematic for two reasons. Firstly, $\mathfrak X$ might be finite and thus such a constant might not exists. Secondly, and more seriously, such a choice would heavily depend on $\mathfrak X$ and would not be natural from a category-theoretic standpoint that demands are purely finitistic and syntactic treatment.}

 \item  If $E\equiv \forall \xi \,F(\xi)$ then
      $$\Gamma_{i+1}\;=\; \Gamma_i',F(\bar{\delta}),\Gamma_i'',\neg A_i$$
      where  $\bar{\delta}$ is any ordinal constant.

      \item If  $E\equiv (\exists \xi<\bar{\alpha})\,F(\xi)$. 
then
      $$\Gamma_{i+1}\;=\; \Gamma_i',\bar{\delta}<\bar{\alpha}\,\wedge\,F(\bar{\delta}),\Gamma_i'',\neg A_i,E$$
      where the ordinal constant $\bar{\delta}$ is picked as follows:  If there occurs an ordinal constant $\bar{\gamma}$ in  $\Gamma_0;\ldots;\Gamma_i$ such that $\bar{\gamma}<\bar{\alpha}$ is true and $F(\bar{\gamma})$ does not
occur in $\Gamma_0;\ldots;\Gamma_i$, let $\bar{\delta}$ be the least such. 
If, however, for all ordinal constants  $\bar{\gamma}$ in  $\Gamma_0;\ldots;\Gamma_i$ with $\bar{\gamma}<\bar{\alpha}$  true the formula $F(\bar{\gamma})$ already
occurs in $\Gamma_0;\ldots;\Gamma_i$, then $\bar{\delta}$ can be any ordinal constant. 

 \item  If $E\equiv (\forall \xi<\bar{\alpha}) \,F(\xi)$ then
      $$\Gamma_{i+1}\;=\; \Gamma_i',\bar{\delta}<\bar{\alpha}\to F(\bar{\delta}),\Gamma_i'',\neg A_i$$
      where  $\bar{\delta}$ is any ordinal constant.

 \item If $E\equiv P^{\bar{\alpha}}_A(s)$. 
then 
      $$\Gamma_{i+1}\;=\; \Gamma_i',A( P^{<\bar{\alpha}}_A,s), \Gamma_i'',\neg A_i.$$
 \item Let $E\equiv \neg P^{\bar{\alpha}}_A(s)$. 
Then 
      $$\Gamma_{i+1}\;=\; \Gamma_i',\neg A( P^{<\bar{\alpha}}_A,s), \Gamma_i'',\neg A_i.$$

      \end{enumerate}
\end{itemize}
 }\end{deff}

  The set of $\PAO$-deduction chains forms a tree ${\mathbb T}_{\mathfrak X}^{^{\bar{Q}}}$ labeled with
strings of sequents. We will now consider two cases.
\begin{itemize}
\item  {\bf Case I:} There is a wellordering $\mathfrak X$ such that  ${\mathbb T}_{\mathfrak X}^{^{\bar{Q}}}$   is {\bf not} well-founded.
 \item[]  Then ${\mathbb T}_{\mathfrak X}^{^{\bar{Q}}}$ contains an infinite path $\mathbb P$.
 \\[2ex]
\item  {\bf Case II}: All ${\mathbb T}_{\mathfrak X}^{^{\bar{Q}}}$  are well-founded.
\end{itemize}

 \paragraph{Case I:}
Let $\mathbb P$ be an infinite path through  ${\mathbb T}_{\mathfrak X}^{^{\bar{Q}}}$. Let $\mathsf{Ord}_{\mathbb P}$ be the set of ordinal constants that occur in sentences of $\mathbb P$.
The language $L_{\mathbb P}$ is the restriction of $\LQ$ that has only the ordinal constants in $\mathsf{Ord}_{\mathbb P}$. We now define a structure $\mathfrak{M}_{\mathbb P}$ for this language.
The number-theoretic part of $\mathfrak{M}_{\mathbb P}$ is just the standard model of the naturals plus the interpretation of $Q$ as $\bar{Q}$ while the ordinal part has as its universe $\mathsf{Ord}_{\mathbb P}$ with the ordering:
\begin{eqnarray*}\bar{\alpha}<_{_{\mathfrak{M}_{\mathbb P}}}\bar{\beta} &\mbox{ iff }& \bar{\alpha}<\bar{\beta}\mbox{ is a true literal}.\end{eqnarray*}
Note that $<_{_{\mathfrak{M}_{\mathbb P}}}$ is a wellordering because $\mathfrak X$ is.

Also the predicates $P_A$ have to be furnished with an interpretation:
\begin{eqnarray*} P_A^{_{\mathfrak{M}_{\mathbb P}}}(n,\bar{\beta}) &\mbox{ iff }& \mbox{ $\neg P_A(\bar{\beta},\bar{n})$ is a formula of a sequent occurring in $\mathbb P$.}\end{eqnarray*}
The aim is to show that $\mathfrak{M}_{\mathbb P}$ models $\mathbf{PA}_{\mathsf{Ord}_{\mathbb P}}$.

\begin{deff}{\em For a sentence $F$ of $L_{\mathbb P}$ we write $F\in \mathbb{P}$ if $F$ occurs in a sequent that belongs to $\mathbb P$.}\end{deff}

\begin{lem} \label{Zerlegung}\begin{enumerate}
\item $\mathbb P$ does not contain any true literals.

\item If $\mathbb P$ contains $E_0\lor E_1$ then $\mathbb P$  contains $E_0$ and $E_1$.
\item If $\mathbb P$ contains $E_0\land E_1$ then $\bbM{P}$ contains $E_0$ or $E_1$.
\item If $\mathbb P$ contains $\exists xF(x)$ then $\bbM{P}$ contains $F(\bar{n})$ for all $n$.
\item If $\mathbb P$ contains $\forall xF(x)$ then $\bbM{P}$ contains $F(\bar{n})$ for some  $n$.
\item If $\mathbb P$ contains $\exists \xi F(\xi)$ then $\bbM{P}$ contains $F(\bar{\beta})$ for all  $\bar{\beta}\in \mathsf{Ord}_{\mathbb P}$.
\item If $\mathbb P$ contains $\forall \xi F(\xi)$ then $\bbM{P}$ contains $F(\bar{\beta})$ for some  $\bar{\beta}\in \mathsf{Ord}_{\mathbb P}$.
\item If $\mathbb P$ contains $(\exists \xi <\bar{\alpha})F(\xi)$ then $\bbM{P}$ contains $F(\bar{\beta})$ for all  $\bar{\beta}\in \mathsf{Ord}_{\mathbb P}$ such that $\bar{\beta}<\bar{\alpha}$ holds true.
\item If $\mathbb P$ contains $(\forall \xi<\bar{\alpha}) F(\xi)$ then $\bbM{P}$ contains $F(\bar{\beta})$ for some  $\bar{\beta}\in \mathsf{Ord}_{\mathbb P}$.
\item If $\bbM{P}$ contains $\neg P_A^{\bar{\alpha}}(s)$ then  $\bbM{P}$ contains $\neg A(P^{<\bar{\alpha}},s)$.
\item If $\bbM{P}$ contains $ P_A^{\bar{\alpha}}(s)$ then  $\bbM{P}$ contains $ A(P^{<\bar{\alpha}},s)$.
\end{enumerate}

\end{lem}
\prf The proof is routine. For an existential formula note that such a formula will become the redex of a formula infinitely many times as it never vanishes.Thus,  for every formula $\exists xF(x)$ in $\mathbb P$, the substitution  instance $F(\bar{n})$ will be in $\mathbb P$ for every $n$. Likewise, for an existential formula $\exists \xi F(\xi)$ in $\mathbb P$,  the substitution  instance $F(\bar{\zeta})$ will be in $\mathbb P$ for every $\bar{\zeta}$ in  $\mathsf{Ord}_{\mathbb P}$; although not necessarily for all ordinal constants $\bar{\eta}$. Similar consideration apply to bounded existential ordinal quantifiers. \qed

\begin{deff}{\em For the next Lemma (or rather its proof) it is convenient to introduce a measure for the complexity of  a $L_{\mathbb P}$-sentences.

An ordinal constant $\bar{\alpha}$ is of the form $\alpha_u$ for a unique $u\in X$. Define its rank, $rk(\bar{\alpha})$, to be the order-type of the initial segment  of $<_X$ determined by $u$.

Let $u_0$ be the least element of $X$ with regard to $<_X$.
\begin{enumerate} 
\item  $|F|=0$ if $F$ is a literal not containing ordinal constants.
\item $|\bar{\alpha}<\bar{\beta}|= |\neg \bar{\alpha}<\bar{\beta}|=\max(rk(\bar{\alpha}),rk(\bar{\beta}))\cdot\omega$. 

\item $|\bar{\alpha}=\bar{\beta}|=|\neg \bar{\alpha}=\bar{\beta}| =\max(rk(\bar{\alpha}),rk(\bar{\beta}))\cdot\omega$. 



\item $|P_A(s,\bar{\alpha})|=|\neg P_A(s,\bar{\alpha})| = rk(\bar{\alpha})\cdot \omega+\omega$.

\item  $|B\wedge C|=|B\vee C|=\max(|B|,|C|)+1$.

\item  $|\exists xF(x)|=|\forall xF(x)|=|F(\bar{0})|+1$.

\item $|\exists\xi F(\xi)|=|\forall \xi F(\xi)|=\max(\tau\cdot\omega, |F(\alpha_{u_0})|+1)$ where $\tau$ is the order-type of $\mathfrak X$.

\item $|(\exists\xi<\bar{\alpha}) F(\xi)|=|(\forall \xi <\bar{\alpha})F(\xi)|=\max(rk(\bar{\alpha})\cdot\omega, |F(\alpha_{u_0})|+1)$.
 
\end{enumerate}

}\end{deff}

\begin{lem}\label{neg} If $F\in \mathbb{P}$ then $\mathfrak{M}_{\mathbb P}\models \neg F$. \end{lem}
\prf The proof proceeds by induction on $|F|$. The claim is obvious for literals as $\mathbb P$ does not contain true literals. 

If $F$ is of the form  $\neg P_A(s,\bar{\alpha})$ then, by definition of  $\mathfrak{M}_{\mathbb P}$, we have  $\mathfrak{M}_{\mathbb P}\models P_A(s,\bar{\alpha})$, whence  $\mathfrak{M}_{\mathbb P}\models \neg F$.

Let be $F$ is of the form  $P_A(s,\bar{\alpha})$.  Then, by Lemma \ref{Zerlegung}, $ A(P^{<\bar{\alpha}},s)\in \mathbb{P}$. One also has that  $|A(P^{<\bar{\alpha}},s)|<|P_A(s,\bar{\alpha})|$, and hence, by the induction hypothesis,
$(*)\;\;\mathfrak{M}_{\mathbb P}\models \neg A(P^{<\bar{\alpha}},s)$. If  $\mathfrak{M}_{\mathbb P}\models P_A(s,\bar{\alpha})$ were to hold, this would mean that  $\neg P_A(s,\bar{\alpha})\in \mathbb{P}$ and thus
$\neg A(P^{<\bar{\alpha}},s)\in \mathbb{P}$ by Lemma \ref{Zerlegung}, yielding $\mathfrak{M}_{\mathbb P}\models  A(P^{<\bar{\alpha}},s)$, again by invoking the induction hypothesis, which contradicts $(*)$. 
As a result, $\mathfrak{M}_{\mathbb P}\models P_A(s,\bar{\alpha})$.

If $F$ is of the form $\exists\xi G(\xi)$ then $G(\bar{\beta})\in \mathbb{P}$ for all  $\bar{\beta}\in \mathsf{Ord}_{\mathbb P}$, and hence, by the induction hypothesis, $\mathfrak{M}_{\mathbb P}\models \neg G(\bar{\beta})$ holds for
all $\bar{\beta}\in \mathsf{Ord}_{\mathbb P}$, thence $\mathfrak{M}_{\mathbb P}\models \neg F$.

If $F$ is of the form $\forall\xi G(\xi)$ then $G(\bar{\beta})\in \mathbb{P}$ for some   $\bar{\beta}\in \mathsf{Ord}_{\mathbb P}$, and hence, by the induction hypothesis, $\mathfrak{M}_{\mathbb P}\models \neg G(\bar{\beta})$, thence $\mathfrak{M}_{\mathbb P}\models \neg F$.

The remaining cases are similar.
\qed

\begin{cor}\label{piein} \begin{itemize}
\item[(i)]  $\mathfrak{M}_{\mathbb P}\models\forall\alpha \forall x\,[A(P_A^{<\alpha},x)\leftrightarrow P_A^{\alpha}(x)]$.
\item[(ii)]  $\mathfrak{M}_{\mathbb P}\models \forall x\,[A(P_A,x)\leftrightarrow P_A(x)]$.
\item[(iii)] If $U\subseteq\mathbb{N}$ is a set such that $\forall x\, [A(U,x)\rightarrow U(x)]$, then $\{n\in\mathbb{N}\mid  \mathfrak{M}_{\mathbb P}\models P_A(\bar{n})\}\subseteq U$.
\end{itemize}
\end{cor} 
\prf (i): The formula $\forall\alpha \forall x\,[A(P_A^{<\alpha},x)\rightarrow P_A^{\alpha}(x)]$ is a stage axiom, and hence its negation $\neg \forall\alpha \forall x\,[A(P_A^{<\alpha},x)\rightarrow P_A^{\alpha}(x)]$ is in $\mathbb{P}$, thus
$\mathfrak{M}_{\mathbb P}\models\forall\alpha \forall x\,[A(P_A^{<\alpha},x)\to P_A^{\alpha}(x)]$. 

Conversely, if $\mathfrak{M}_{\mathbb P}\models P_A^{\bar{\beta}}(s)$ then  $\neg P_A^{\bar{\beta}}(s) \in\mathbb{P}$, thus $\neg A(P^{\bar{\alpha}},s)\in \mathbb{P}$, yielding 
$\mathfrak{M}_{\mathbb P}\models A(P^{\bar{\alpha}},s)$ by Lemma \ref{neg}. 
\\[1ex] \indent
(ii): Note that the forward part of the equivalence is a reflection axiom.  So the formula $\neg \forall\alpha \forall x\,[A(P_A,x)\to P_A(x)]$ is in $\mathbb{P}$, and therefore  $\mathfrak{M}_{\mathbb P}\models\forall\alpha \forall x\,[A(P_A,x)\to P_A(x)]$ by Lemma \ref{neg}. If   $\mathfrak{M}_{\mathbb P}\models P_A(s)$, then  $\mathfrak{M}_{\mathbb P}\models P_A^{\bar{\beta}}(s)$ for some $\bar{\beta}$. The latter entails  
$\mathfrak{M}_{\mathbb P}\models A(P_A^{<\bar{\beta}},s)$ by (i), which on account of positivity implies $\mathfrak{M}_{\mathbb P}\models A(P_A,s)$.
\\[1ex] \indent
(iii): Assume $\forall x\, [A(U,x)\rightarrow U(x)]$.  By induction over $\mathsf{Ord}_{\mathbb P}$  we prove $$\{n\in\mathbb{N}\mid  \mathfrak{M}_{\mathbb P}\models P_A^{\bar{\alpha}}(\bar{n})\}\subseteq U,$$  from which (iii) follows.
By induction hypothesis,  $$\{ n\in\mathbb{N}\mid \mathfrak{M}_{\mathbb P}\models P_A^{\bar{\beta}}(\bar{n}\}\subseteq U$$ holds for all $\bar{\beta}$ such that $\mathfrak{M}_{\mathbb P}\models \bar{\beta}<\bar{\alpha}$.
Thus, $$\{ n\in\mathbb{N}\mid  \mathfrak{M}_{\mathbb P}\models A( P_A^{<\bar{\alpha}},\bar{n})\}\subseteq \{n\in\mathbb{N}\mid A(U,n)\}$$  by positivity, and hence, by (i), 
$\{n\in\mathbb{N}\mid  \mathfrak{M}_{\mathbb P}\models P_A^{\bar{\alpha}}(\bar{n})\}\subseteq U$.
\qed

\begin{cor}\label{Pieinsca}{\em  ($\mathbf{ACA}_0$)} The class $\{n\in \mathbb{N}\mid H(n,\bar{Q})\}$ with $H(x,U)$ being a $\Pi^1_1$-formula and all free variables exhibited, is first-order definable over $\mathfrak{M}_{\mathbb P}$. Thus it is a set.
\end{cor}
\prf It is well-known\footnote{See for instance \cite[III.3.2]{hinman} or \cite[V.1.8]{SOSA}.} (and provable in $\mathbf{ACA}_0$) that every class $\{n\in \mathbb{N}\mid H(n,\bar{Q})\}$  is many-one reducible to one of the form 
\begin{eqnarray}  & & \{n\in \mathbb{N}\mid \forall X( \forall u[A(X,u)\to u\in X] \to n\in X)\}.\end{eqnarray}
The latter class, however, is the same as $\{n\in\mathbb{N}\mid \mathfrak{M}_{\mathbb P}\models P_A(\bar{n})\}$ by 
Corollary \ref{piein}. \qed

From the foregoing we conclude that $\Pi^1_1$-comprehension obtains (on the basis of $\mathbf{ACA}_0$)  if for all sets of naturals $Q$ there exists a wellordering $\mathfrak X$ such that the the search tree  ${\mathbb T}_{\mathfrak X}^{^{\bar{Q}}}$   is 
ill-founded. Thus we want to exclude the possibility that for some $\bar{Q}$ all of the search trees ${\mathbb T}_{\mathfrak X}^{^{\bar{Q}}}$  for wellorderings $\mathfrak X$  are well-founded. 
To this we turn next.

\paragraph{  Case II: All ${\mathbb T}_{\mathfrak X}^{^{\bar{Q}}}$  are well-founded.} The strings of sequents that make up the nodes of ${\mathbb T}_{\mathfrak X}^{^{\bar{Q}}}$ can be linearly ordered by the Kleene-Brouwer 
tree ordering, which is a wellordering on account of the well-foundedness of the tree. 
Let $\mathsf{KB}^{^{\bar{Q}}}_{ \mathfrak{X}}$ be the Kleene-Brouwer ordering of ${\mathbb T}_{\mathfrak X}^{^{\bar{Q}}}$ (see e.g. \cite[V.1.2]{SOSA} for a definition).
Owing to the careful design of  ${\mathbb T}_{\mathfrak X}^{^{\bar{Q}}}$, the map $\mathfrak{X}\mapsto ({\mathbb T}_{\mathfrak X}^{^{\bar{Q}}},\mathsf{KB}^{^{\bar{Q}}}_{ \mathfrak{X}})$
gives rise to an endofunctor $F_{\!_{\bar{Q}}}$ of the category of wellorderings. To see this, let $\mathfrak Y$ be another wellordering and $f\co\mathfrak{X}\to\mathfrak{Y}$ be a an order preserving function. 
Each element of ${\mathbb T}_{\mathfrak X}^{^{\bar{Q}}}$ can be viewed as a notation $S(\alpha_{u_0},\ldots,\alpha_{u_{n-1}})$ where $\alpha_{u_0},\ldots,\alpha_{u_{n-1}}$ are the ordinal constants occurring in it with
$u_0,\ldots, u_{n-1}\in X$ and $u_0<_X\ldots<_X u_{n-1}$. Now define $F_{\!_{\bar{Q}}}(f)$ by 
$$F_{\!_{\bar{Q}}}(f) (S(\alpha_{u_0},\ldots,\alpha_{u_{n-1}})) \;:=\; S(\alpha_{f(u_0)},\ldots,\alpha_{f(u_{n-1})}).$$
Then $$F_{\!_{\bar{Q}}}(f): ({\mathbb T}_{\mathfrak X}^{^{\bar{Q}}},\mathsf{KB}^{^{\bar{Q}}}_{ \mathfrak{X}})\to ({\mathbb T}_{\mathfrak Y}^{^{\bar{Q}}},\mathsf{KB}^{^{\bar{Q}}}_{ \mathfrak{Y}}).$$
We shall denote the denotation system associated with this dilator by ${\mathcal D}_{\!{\bar{Q}}}$.
More concretely, ${\mathcal D}_{\!{\bar{Q}}}$ can be visualized as consisting of the denotations
   $$(S(x_0,\ldots,x_{n-1});u_0,\ldots,u_{n-1};\mathfrak{X})$$
where  $S(\alpha_{u_0},\ldots,\alpha_{u_{n-1}})\in {\mathbb T}_{\mathfrak X}^{^{\bar{Q}}}$.

For notational reasons, let's drop the subscript $\bar{Q}$, writing $\mathcal{D}:={\mathcal D}_{\!{\bar{Q}}}$. 
Next, $\mathcal{D}$ is wedded to the Bachmann construction as in Definition \ref{Bachdil}, giving rise to the ordinal representation system 
$\BHSD$.
The aim is now to show that the well-foundedness of $\BHSD$ leads to a contradiction.
Let $\mathfrak{I}$ be the subordering of $\BHSD\cap\Omega$ of $\BHSD$ obtained by allowing only elements $\alpha$ with $\alpha<\Omega$.
By assumption, the search tree ${\mathbb T}_{\mathfrak I}^{^{\bar{Q}}}$ is well-founded. 
A first step consists in realizing that  ${\mathbb T}_{\mathfrak I}^{^{\bar{Q}}}$ is actually contained in $\BHSD$. Moreover, this order-morphism sends elements 
 of ${\mathbb T}_{\mathfrak I}^{^{\bar{Q}}}$ to $\varepsilon$-numbers above $\Omega$ in $\BHSD$. 

\begin{lem}\label{rein} There is an order-morphism that maps ${\mathbb T}_{\mathfrak I}^{^{\bar{Q}}}$, equipped with the Kleene-Brouwer ordering,  into $\BHSD$.
Moreover, this order-morphism sends elements 
 of ${\mathbb T}_{\mathfrak I}^{^{\bar{Q}}}$ to $\varepsilon$-numbers above $\Omega$ in $\BHSD$. 
\end{lem}
\prf An element of ${\mathbb T}_{\mathfrak I}^{^{\bar{Q}}}$ is a string of sequents, $S(\alpha_0,\ldots,\alpha_{n-1})$,  whose ordinal constants  $\alpha_0<\ldots<\alpha_{n-1}$ are from ${\mathfrak I}$.
Letting $c:=S(x_0,\ldots,x_{n-1})$ and $\sigma:= D(c;0,\ldots,n-1;n)$, define the desired map by
$$S(\alpha_0,\ldots,\alpha_{n-1})\mapsto \mathfrak{E}^{\sigma}_{\alpha_0,\ldots,\alpha_{n-1}}.$$
By definition of $\mathcal D$ and the ordering of $\BHSD$ it is (quite) clear that this mapping preserves the order.
\qed

An important step consists in viewing ${\mathbb T}_{\mathfrak I}^{^{\bar{Q}}}$ as the skeleton of a proof 
in an infinite proof system that enjoys partial cut elimination, for this will mean that one obtains a cut-free proof of the empty sequent, which is impossible. 
The infinitary system we have in mind is a variant of $\PAIs$. To be precise,  Definition \ref{Axiome} only lists the axioms and subsection \ref{DedC} only defines deductions chains in $\PAIs$. To render it a proper proof system, we have to add  inference  rules.

\begin{deff}\label{Rules}{\em  $\PAI$ has the same axioms as $\PAIs$ except for the stage and reflection axioms, which will be turned into rules.

 The inference rules of $\PAI$ are:\footnote{There is a slightly annoying aspect to this proof system. $\PAI$ is supposed to be a proper deduction calculus, which entails that we have to pay attention 
 to structural rules, too. These boring rules can be moved under the carpet by conceiving of sequents as finite sets of formulas. Thus,  in what follows,  sequents in the sense of 
 $\PAI$ will be understood as finite sets of formulas. Notations like $\Gamma,A$ and $\Gamma,\Lambda$ will consequently have to be deciphered as $\Gamma\cup\{A\}$ and $\Gamma\cup\Lambda$, respectively.}
$$\begin{array}{ll}
(\wedge)&  {\frac{\DI \Gamma, E_0 \;\;\;\Gamma, E_1}{{\DI \Gamma, E_0
\wedge E_1^{\phantom{i}}}}^{}}\\[0.6cm]

(\vee)&{\frac{\DI \Gamma, E_{i}}{\DI\Gamma, E_{0} \vee E_{1}^{\phantom{i}}}}\;\;
\mbox{ if } i=0 \mbox{ or } i=1\\[0.6cm]

(\forall_{\!_{\mathbb N}})& {\frac{\DI\Gamma,\, F(\bar{n})\;\;\mbox{ for all $n\in\mathbb{N}$ }}
{\DI\Gamma,( \forall u\in\mathbb{N})
\, F(u)^{\phantom{i}}}}\\[0.6cm]
(\exists_{_{\mathbb N}})& {\frac{\DI\Gamma, \, F(s)}{\DI\Gamma, (\exists u\in\mathbb{N})
\,F(u)^{\phantom{i}}}} \;\;\mbox{  $s$ number term}
\\[0.6cm]

(b\forall)& {\frac{\DI\Gamma,  F(\bar{\beta})\;\;\mbox{ for all $\bar{\beta}<\bar{\alpha}$}}
 {\DI\Gamma, (\forall \xi<\bar{\alpha}) F(\xi)^{\phantom{i}}}}\\[0.6cm]

(b\exists)& {\frac{\DI\Gamma,  F(\bar{\beta})}{\DI\Gamma,
(\exists \xi<\bar{\alpha})F(\xi)^{\phantom{i}}}}\;\; \mbox{ if } \bar{\beta}<\bar{\alpha}\\[0.6cm]

(\forall)& {\frac{\DI\Gamma,\, F(\bar{\alpha})\;\;\mbox{ for all $\bar{\alpha}$ }}
{\DI\Gamma, \forall \beta
\, F(\beta)^{\phantom{i}}}}\\[0.6cm]
(\exists)& {\frac{\DI\Gamma, \, F(\bar{\alpha})}{\DI\Gamma, \exists \beta
\,F(\beta)^{\phantom{i}}}}
\\[0.6cm]
(P_A)& {\frac{\DI\Gamma,  A(P_A^{<\bar{\alpha}},s)}{\DI\Gamma,
P_A(\bar{\alpha},s)^{\phantom{F}}}}\\[0.6cm]

(\neg P_A)& {\frac{\DI\Gamma,  \neg A(P_A^{<\bar{\alpha}},s)}{\DI\Gamma,
\neg P_A(\bar{\alpha},s)^{\phantom{F}}}}\\[0.6cm]

(\mathrm{Ref})&{\frac{\DI\Gamma, A(P_A,s)}{\DI\Gamma, \,\exists \xi\, P_A(\xi,s)^{\phantom{i}}}}\;\;\\[0.6cm]

\Cut& {\frac{\DI\Gamma, A \;\;\;\;\;\; \Gamma,
\neg A}{\DI\Gamma^{\phantom{i}}}}

\end{array}$$ 
where in the above rules $s$ is always a number term.
}\end{deff}

The first observation worthy of note  is that any sequent provable in ${\mathbb T}_{\mathfrak I}^{^{\bar{Q}}}$ is provable in $\PAI$ in a very uniform way.
${\mathbb T}_{\mathfrak I}^{^{\bar{Q}}}$ consists of deduction chains $\Gamma_0;\ldots; \Gamma_{r-1}$, and we say that $\Gamma_{r-1}$ is {\em provable in ${\mathbb T}_{\mathfrak I}^{^{\bar{Q}}}$}. More precisely, we write
$$\TsQ{S(\alpha_0,\ldots,\alpha_{n-1})}{}{\Gamma}$$
if there exists a deduction chain $\Gamma_0;\ldots; \Gamma_{r-1}$ such that $S(\alpha_0,\ldots,\alpha_{n-1})$ is the notation for this chain and $\Gamma$ is $\Gamma_{r-1}$.

It has always been a challenge of impredicative ordinal-theoretic proof theory to control infinitary derivations. They have to be of a very uniform kind so as to be able to prove  (partial) cut elimination for such derivations.  The problem arises for first time when analyzing systems such as $\mathbf{ID}_1$, $\mathbf{KP}$ and bar induction whose proof-theoretic ordinal 
is the Bachmann-Howard ordinal. There are several technical approaches from the 1980s,  e.g.  \cite{j80,pohlers81,j82,bu,pohlers}, that  could be enlisted to analyze $\PAI$, too.
These days, though,  Buchholz'  technique of operator controlled from \cite{bu:92} is most widely used and  we will also refer to it when defining the uniform derivations of $\PAI$.
However, we will not go into much detail as these techniques are standard by now.\footnote{E.g. in  the textbook \cite{pohlers09}.}  
For our purposes, an {\em operator} $\mathcal H$ is just a  monotone, inclusive and
 idempotent function from  subsets of  $\BHSD$ to subsets of  $\BHSD$ such that for all $X$, $\mathcal{H}(X)$ contains $0$ and $\Omega$ and the following holds:
 \begin{itemize}
 \item If $\alpha$ has Cantor normal form $\omega^{\alpha_1}+\ldots+\omega^{\alpha_n}$ then $\alpha\in\mathcal{H}(X)\;\Leftrightarrow \alpha_1,\ldots,\alpha_n\in \mathcal{H}(X)$;
 \item If  $\alpha_0,\ldots,\alpha_{n-1}\in X$, $\alpha_0<\cdots<\alpha_{n-1}<\Omega$, and 
$\sigma=(c;0,\ldots,n-1;n)\in \mathcal{D}$, then
$\EDD^{\sigma}_{\alpha_0,\ldots,\alpha_{n-1}}\in\mathcal{H}( X).$ 
\end{itemize}
We shall write $$ \provx{{\mathcal H}}{\delta}{\rho}{\Gamma}$$ to convey that the sequent $\Gamma$ is deducible in $\PAI$ controlled by the operator $\mathcal H$ with length $\delta$ and cut rank $\rho$.

For an operator $\mathcal H$ and  finite subset $U$ of $\BHSD$, $\mathcal{H}[U]$ stands for the operator
with $\mathcal{H}[U](V) := \mathcal{H}(U\cup V)$.
One then obtains the following interpretation result with the length of the derivation determined by the order-morphism of Lemma \ref{rein}. 

\begin{thm}\label{inter} For every operator $\mathcal H$, \begin{eqnarray*} \TsQ{S(\alpha_0,\ldots,\alpha_{n-1})}{}{\Gamma} &\Rightarrow &  \provx{{\mathcal H}[\{\alpha_0,\ldots,\alpha_{n-1}\}]}{\mathfrak{E}^{\sigma}_{\alpha_0,\ldots,\alpha_{n-1}}}{\Omega+k}{\Gamma}\end{eqnarray*} for some $k<\omega$,
where $c:=S(x_0,\ldots,x_{n-1})$ and $\sigma:= D(c;0,\ldots,n-1;n)$.

\end{thm}

With $\mathcal H$ controlled derivations one obtains partial cut elimination in that all provable sequents consisting entirely of $\Sigma$-formulas, i.e. formulas without unbounded universal ordinal quantifiers, have a cut-free proof.
The bottom node of the tree  ${\mathbb T}_{\mathfrak I}^{^{\bar{Q}}}$ is its maximum element in the Kleene-Brouwer ordering. Denoting it by $S_0$, we have $\TsQ{S_0}{}{\emptyset}$
with $\emptyset$ standing for the empty sequent (without any formulas). Owing to the interpretation theorem \ref{inter}, we then have 
\begin{eqnarray}\label{zero}&&\provx{{\mathcal H}}{\mathfrak{E}^{\sigma_0}}{\Omega+3}{\emptyset}\end{eqnarray}
where $c:=S_0$ and $\sigma_0:= D(c)$. 
Now, it is possible to turn the derivation of (\ref{zero}) into a cut-free one. Technically, however, the cut elimination procedure requires a hierarchy of ever stronger operators ${\mathcal H}^*_{\xi}$  with $\xi\in\BHSD$ that, with growing ordinal index,  accommodate more and more ordinals engendered by the collapsing 
function $\varthetad$, which is to say that the sets ${\mathcal H}^*_{\xi}(X)$ enjoy ever stronger closure properties with regard to $\varthetad$ as $\xi$ increases.
As this is a standard technique in ordinal analysis (see \cite{bu:92,pohlers09,Ra12}) at the lowest level of impredicativity, we will not go into details.
The upshot is that by applying cut elimination to the derivation in (\ref{zero}) we arrive at a cut-free $\PAI$-derivation of the empty sequent
$$\provx{{\mathcal H}^*_{\rho}}{\tau}{0}{\emptyset}$$
for  specific $\rho,\tau\in \BHSD$  with $\rho<\Omega$.  
But this is absurd for a cut-free proof as (Cut) is the only inference that could yield the empty sequent as a consequence; hence case II has been ruled out, thus it's always case I,
concluding this direction of the main theorem's proof.

  \subsection{A glimpse of Anton Freund's work}

In his thesis and a series of papers \cite{freund-thesis, freund0, freund1,freund2},  Anton Freund presents a much broader perspective on the connection between 
$\Pi^1_1$-Comprehension  and various Bachmann-Howard principles. In particular, he developed a nice categorical approach to the support operation as a natural transformation, thereby introducing the notion of a
prae-dilator.

\begin{deff} {\em The functor $$[\,\cdot\,]^{<\omega}\co \LIO \to \mathbf{Set}$$ is defined by
letting  \begin{eqnarray*} [\mathfrak{X}]^{<\omega} &= & \mbox{ the set of finite subsets of $X$}
\\[1ex] [f]^{<\omega}\co\; [\mathfrak{X}]^{<\omega} &\to& [\mathfrak{Y}]^{<\omega}
\\[1ex]
[f]^{<\omega}(\{x_1,\ldots,x_n\}) &=& \{f(x_1),\ldots,f(x_n)\}
\end{eqnarray*} whenever $f\co \mathfrak{X}\to\mathfrak{Y}$ and $x_1,\ldots,x_n\in X$. 
 
An endofunctor $T\co \LIO \to \LIO$ is called a {\em prae-dilator} if there exists a natural transformation
$$\supp\co T\Rightarrow [\,\cdot\,]^{<\omega}$$ i.e.,  
 the square 
\[
\xymatrix{
T(\mathfrak{X}) \ar[d]_{T(f)}  \ar[r]^{\supp_{\mathfrak X}} & [\mathfrak{X}]^{<\omega} \ar[d]^{[f]^{<\omega}} \\
T(\mathfrak{Y}) \ar[r]_{\supp_{\mathfrak Y}} & [\mathfrak{Y}]^{<\omega}}
 \]
commutes for every $f\co \mathfrak{X}\to \mathfrak{Y}$ in $\LIO$,
and for $\sigma\in T(\mathfrak{X})$, 
$$\sigma \in\ran( T(\iota_{\sigma})),$$
where $\iota_{\sigma}\co \supp_{\mathfrak X}(\sigma)\hookrightarrow \mathfrak{X}$ denotes the obvious inclusion map.
}\end{deff}
It's not the purpose of this paper to present more details of Freund's work, however,  we strongly recommend looking at the original sources  \cite{freund-thesis, freund0, freund1,freund2}.
To finish this subsection, we quote the following characterization theorem from \cite{freund-thesis, freund1,freund2}.
\begin{thm} The following are equivalent over $\mathbf{RCA}_0$:
 
 \begin{itemize}
 \item[(i)] The principle of $\Pi^1_1$-comprehension.
 \item[(ii)] The abstract Bachmann-Howard principle: {\em Every dilator has a well-founded Bachmann-Howard fixed point}.
 
\item[(iii)] The computable Bachmann-Howard principle: {\em For every dilator $T$ the linear order $\vartheta(T)$
is well-founded.}
\end{itemize}
 \end{thm}

 \section{There are much stronger constructions than Bachmann's}
In the list of strong theories for which an ordinal analysis has been achieved, the theory of non-iterated inductive arithmetic definitions (embodied in $\mathbf{PA}_{\Omega}$)  is just the starting point. Likewise, Bachmann's ordinal representation system is just the first
 in a long line of ever bigger ones. It is thus natural to ponder what happens when one marries the stronger representation systems with a dilator in the same vein as \ref{Bachdil}.
As some $\WOP(f)$ principles turned out  to be equivalent to the existence of countable coded $\omega$-models, one is led to think that there is a stronger notion of model 
 that pertains to principles such as $\WOPP(F)$.
  After $\omega$-models come $\beta$-models and we take  interest in statements of the form  {\em ``every set belongs to a countably coded $\beta$-model of $T$'' for various theories $T$}.
 \comment{ and the theory $\Pi^1_1\mbox{-}\mathbf{CA}$ has a characterization
  in terms of countable coded $\beta$-models, namely
   via the statement ``every set belongs to a countably coded $\beta$-model''. }

\begin{deff}\label{models} {\em Every $\omega$-model ${\mathfrak M}$ of a theory $T$
in the language of second order arithmetic is isomorphic to a structure  $${\mathfrak A}=\langle \omega;{\mathfrak X};0,+,\times,\ldots\rangle$$ where
${\mathfrak X}\subseteq {\mathcal P}(\omega)$.

 ${\mathfrak A}$ is a {\em $\beta$-model} if the concept
of well ordering is absolute with respect to ${\mathfrak A}$,
i.e. for all $X\in {\mathfrak X}$,
 \begin{eqnarray*} {\mathfrak A}\models \WO(<_{X}) &\mbox{ iff }&
  \mbox{ $<_X$ is a well ordering},\end{eqnarray*}
  where $n<_Xm :\Rightarrow 2^n3^m\in X$.
 } \end{deff}
A natural test case is to consider  $\beta$-models of $\Pi^1_1$-comprehension. There already exists a characterization of the statement that every set is contained in an $\omega$-model of 
$\Pi^1_1$-comprehension in \cite{rt} due to  Ian Alexander Thomson and the author of this paper, as well as a characterization of the statement that every set is contained in an $\omega$-model of 
$\Pi^1_1$-comprehension plus bar induction due to  Ian Alexander Thomson \cite{alec}.
This requires relativizing a stronger ordinal representation system than Bachmann's. The construction  gives rise to a dilator, too.

\begin{deff}\label{pirel} {\em 
We do not want to elaborate greatly on the details of the definition of the representation system $\mathsf{OT}(\Omega_{\omega}\cdot{\fX})$ for a wellordering $\fX$,  but like  to give the main ideas.  The crucial extension here is that there is not just one strong ordinal $\Omega$ as in Definition \ref{Bach} but infinitely many
$\Omega_1,\Omega_2,\ldots,\Omega_n,\ldots$  with $\Omega_{\omega}$ denoting their limit. $\fX$ enters the representation system $\mathsf{OT}(\Omega_{\omega}\cdot{\fX})$
as follows: For $u\in \fX$ one has $\Omega_{\omega}\cdot u\in \mathsf{OT}(\Omega_{\omega}\cdot{\fX})$, and if $u<_{\fX}v$ and $\alpha<\Omega_{\omega}$, then
$\Omega_{\omega}+\alpha<\Omega_{\omega}\cdot u$ and  $\Omega_{\omega}\cdot u +\alpha< \Omega_{\omega}\cdot v$.

Moreover, each $\Omega_{n}$ comes equipped with its own (partial)  collapsing function $\vartheta_{n}^{\fX}$, mapping the elements of  $\mathsf{OT}(\Omega_{\omega}\cdot{\fX})$
below $\Omega_{n}$, more precisely
$$\vartheta_{n}^{\fX}: \mathsf{OT}(\Omega_{\omega}\cdot{\fX})\longrightarrow (\Omega_{n-1},\Omega_n)$$
where $ (\Omega_{n-1},\Omega_n)$ denotes the interval of terms strictly inbetween $\Omega_{n-1}$ and $\Omega_n$ with $\Omega_0:=0$.
Furthermore, each  $\Omega_{n}$ has  its own support function $\mathsf{supp}_{\Omega_n}$ and $\vartheta_{n}^{\fX}(\alpha)$ is only defined if 
the ordinals in $\EEON(\alpha)$ are strictly less than $\alpha$, for which we write $\EEON(\alpha)<\alpha$.
$\EEON(\alpha)$ is defined as follows:
\begin{itemize} 
\item[(i)] $\EEON(0)=\emptyset$, $\EEON(\Omega_{\tau})=\emptyset$ for $\tau\in \omega+1$. 
\item[(ii)] $\EEON(\alpha)=\EEON(\alpha_1)\cup\cdots\cup\EEON(\alpha_n)$ if $\alpha=_{CNF}\omega^{\alpha_1}+\cdots+\omega^{\alpha_n}>
\alpha_1$.
\item[(iii)] $\EEON(\Omega_{\omega}\cdot u+ \alpha)=\EEON(\alpha)$ for $\alpha<\Omega_{\omega}$ and $u\in \fX$.
\item[(iv)] For $\beta=\vartheta_m^{\fX}(\gamma)$, $\EEON(\beta)=\emptyset$ if $m< n$ and
$$\EEON(\beta)=\{\gamma\}\cup \EEON(\gamma)$$ if $m \geq n$. 

Finally, the ordering between terms $\vartheta_n^{\fX}(\alpha)$ and $\vartheta_m^{\fX}(\beta)$ is determined as follows:
$$\vartheta_n^{\fX}(\alpha)< \vartheta_m^{\fX}(\beta)\mbox{ iff } n<m\,\vee\,[n=m\,\wedge\,\alpha<\beta]. $$

\end{itemize}
}\end{deff}

\begin{thm} {\em (Rathjen, Thomson \cite{rt}) } Over $\mathbf{RCA}_0$ the following are equivalent:
\begin{itemize}
\item[(i)] Every set is contained in an $\omega$-model of $\Pi^1_1$-comprehension.

\item[(ii)]  $\forall {\fX}\;[\WO(\fX)\rightarrow
\WO(\mathsf{OT}(\Omega_{\omega}\cdot{\fX})]$.
\end{itemize}
\end{thm}
\prf \cite{rt}. \qed

The construction of Definition \ref{pirel} lends itself to a combination  with any dilator $\mathcal D$, following the pattern of Definition \ref{Bachdil}, giving rise to an ordinal representation system
 notated by  $\mathsf{OT}_{\mathcal D}(\Omega_{\omega})$.

\begin{deff}\label{Bachdil2}{\em In $\mathsf{OT}_{\mathcal D}(\Omega_{\omega})$, in place of the objects $\Omega_{\omega}\cdot u$ for $u\in\fX$ one has the following closure property
 of Definition \ref{pirel}:

If  $\alpha_0,\ldots,\alpha_{n-1}\in\mathsf{OT}_{\mathcal D}(\Omega_{\omega})$, $\alpha_0<\cdots<\alpha_{n-1}<\Omega_{\omega}$, and 
$\sigma=(c;0,\ldots,n-1;n)\in \mathcal{D}$, then
$$\EDD^{\sigma}_{\alpha_0,\ldots,\alpha_{n-1}}\in \mathsf{OT}_{\mathcal D}(\Omega_{\omega})\mbox{ and }\Omega_{\omega}<\EDD^{\sigma}_{\alpha_0,\ldots,\alpha_{n-1}}.$$ 
Objects $\EDD^{\sigma}_{\alpha_0,\ldots,\alpha_{n-1}}$ act as $\varepsilon$-numbers above $\Omega_{\omega}$.
The ordering between terms $\EDD^{\sigma}_{\alpha_0,\ldots,\alpha_{n-1}}$ and $\EDD^{\tau}_{\beta_0,\ldots,\beta_{m-1}}$ is defined as in  \ref{pirel}.
Furthermore, the support functions need to be extended to these terms: 
$$\EEON( \EDD^{\sigma}_{\alpha_0,\ldots,\alpha_{n-1}})\;=\;\EEON(\alpha_0)\,\cup\,\ldots\cup\,\EEON(\alpha_{n-1}).$$
Terms  $\vartheta_n^{\mathcal D}(\alpha)$ are generated by the following clause:

If  $\EEON(\alpha)<\alpha$ then $\vartheta_n^{\mathcal D}(\alpha)\in \mathsf{OT}_{\mathcal D}(\Omega_{\omega})$ and $\vartheta_n^{\mathcal D}(\alpha)<\Omega_n$.

Finally, the ordering between terms $\vartheta_n^{\mathcal D}(\alpha)$ and $\vartheta_m^{\mathcal D}(\beta)$ is determined as in \ref{pirel}:
$$\vartheta_n^{\mathcal D}(\alpha)< \vartheta_m^{\mathcal D}(\beta)\mbox{ iff }  n<m\,\vee\,[n=m\,\wedge\,\alpha<\beta].$$

 }\end{deff}
 
 By combining the techniques of \cite{rt,alec} and section 5, one arrives at a characterization of $\beta$-models of $\Pi^1_1$-comprehension.


\begin{thm}\label{BT} Over $\mathbf{RCA}_0$ the following are equivalent:
\begin{itemize}
\item[(i)] Every set is contained in a countably coded $\beta$-model of $\Pi^1_1$-comprehension.
\item[(ii)] For every dilator $\mathcal D$, $\mathsf{OT}_{\mathcal D}(\Omega_{\omega})$ is a well-ordering.
\end{itemize}
\end{thm}
\prf See \cite{Rbeta}.\qed
In all likelihood, Theorem \ref{BT} is not a singularity. For instance, using a representation system with an inaccessible as in the ordinal analysis of $\Delta^1_2$-comprehension plus bar induction (see \cite{jp}), one should arrive at a dilatorial characterization of 
$\beta$-models of $\Delta^1_2$-comprehension. 
Similar characterizations should be obtainable for $\beta$-models of theories axiomatizing recursively Mahlo universes (cf. \cite{r91,rmahlo,Arai2003}) or $\Pi_3$-reflecting universes (cf. \cite{r94b,Arai2004})  and much more, the general conjecture, stated somewhat loosely, reads as follows:

\begin{conj} Ordinal representation systems $\mathsf{OT}$ used in ordinal analyses of strong systems $T$  give rise to functors $$F_{\!_{\mathsf{OT}}}:\mathsf{DIL}\longrightarrow \mathsf{WO}$$
sending dilators to ordinal representation systems. The assertion that $F_{\!_{\mathsf{OT}}}$ sends dilators to wellorderings is equivalent to the statement that every set is contained in a countably coded $\beta$-model of $T$.\end{conj}


\paragraph{Acknowledgements.}  This paper originated from a tutorial given at the Logic Colloquium 2019 in Prague. The author's research was supported by a grant from the John Templeton Foundation
 {\em A new dawn of Intuitionism: Mathematical and Philosophical advances}  (Grant ID 60842).

\end{document}